\documentclass{article}

\usepackage{amssymb}
\usepackage{amsmath}
\usepackage{amscd}
\usepackage{amsthm}
\usepackage{indentfirst}

\begin{document}
\newtheorem{thm}{Theorem}[section]
\newtheorem{cor}[thm]{Corollary}
\newtheorem{propn}[thm]{Proposition}
\newtheorem{defn}[thm]{Definition}
\newtheorem{axe}[thm]{Lemma}
\theoremstyle{definition}
\newtheorem{illus}[thm]{Example}
\newtheorem{remk}[thm]{Remark}
\newtheorem*{tA}{Theorem A}
\newtheorem*{tB}{Theorem B}
\newtheorem*{tC}{Theorem C}
\newtheorem*{tD}{Theorem D}
\newtheorem*{tE}{Theorem E}
\newtheorem*{tF}{Theorem F}

\title{Co-Hopfian Modules }
\author{F. C. Leary  \thanks {St. Bonaventure University (retired), e-mail: chleary@sbu.edu} }
\date{}

\maketitle

\begin{abstract}
If $R$ is a commutative ring with $1$, then a co-Hopfian $R$-module $M$ is necessarily a module of fractions over a ring of fractions. For $R$ commutative Noetherian, we use  Matlis' notion of modules with maximal orders together with localization to obtain results about co-Hopfian submodules of indecomposable injectives and co-Hopfian finitely generated P-primary $R$-modules. With suitable restriction, we detemine co-Hopfian and Hopfian submodules of a co-Hopfian injective module.
\end{abstract}

Heinz Hopf once raised the question: Can a finitely generated group be isomorphic to one of its proper quotients? The question was raised, in a topological context, in 1931 or 1932, depending on your source (\cite[Ch 4 \S 4]{LS} or \cite[\S 2.4]{MKS}, respectively). In \cite{Hopf} he showed that the fundamental group $G$ of a closed two-dimensional orientable surface does not have this property. That is, if $f:G\rightarrow G$ is a surjective group homomorphism, then $f$ must be an isomorphism, and so $G$ is {\em not} isomorphic to a proper quotient of itself. Thus was born the idea of what we now call a {\em Hopfian} group. There is the requisite dual notion of a {\em co-Hopfian} group $G$: if $f:G\rightarrow G$ is an injective group homomorphism, then $f$ must be an isomorphism. Hence, $G$ is not isomorphic to one of its proper subgroups.

According to Goldsmith and Gong \cite{G-G}, R. Baer \cite{Baer} initiated the study of Hopfian and co-Hopfian groups, calling them $Q$-groups and $S$-groups, respectively. The notions Hopfian and co-Hopfian extend in obvious fashion to left (or right) modules over any ring $R$ with $1$ (and to objects in any category).

The Hopfian and co-Hopfian properties are {\em finiteness conditions}. Indeed, in the category of sets, and in the category of vector spaces over a given field $k$, the Hopfian, and the co-Hopfian, objects are precisely the finite sets and the finite dimensisonal $k$-vector spaces, respectively. In general categories, the two concepts are distinct. As an Abelian group, the ring $\mathbb{Z}$ is Hopfian but not co-Hopfian. The  Pr\"{u}fer $p$-group, $\mathbb{Z}(p^{\infty})$, is co-Hopfian but not Hopfian. We will call an $R$-module {\em bi-Hopfian} if it is both Hopfian and co-Hopfian.

We noted in \cite{Me2}, \cite{ Me} the well-known facts that co-Hopfian $R$-modules need not be finite, or even finitely generated, and that finitely generated modules need not be co-Hopfian.  Also, we saw that the class of co-Hopfian $R$-modules fails, in general, to be closed under submodules and quotients.  

Vasconcelos \cite{Wolmer2} and, independently, Strooker \cite {Strooker}, showed that a finitely generated module over a commutative ring is Hopfian. Vasconcelos \cite{Wolmer} also showed that all finitely generated modules over a commutative ring are co-Hopfian if and only if the ring  is 0-dimensional. Hence, all finitely generated modules over a commutative ring are bi-Hopfian if and only if the ring is 0-dimensional. In \cite{Me3}, we showed that over a  commutative Artinian principal ideal ring, a module is Hopfian (co-Hopfian) if and only if it is finitely generated, among other equivalent conditions.

Around 2010, there was renewed interest in Hopfian and co-Hopfian abelian groups because of their connection with the notion of {\em algebraic entropy}, and its dual concept {\em adjoint entropy}, of abelian groups (see, e.g., \cite{DGSZ} and \cite{G-G2}). The papers \cite{G-G}, \cite{Haghany} list some well-known facts about Hopfian and co-Hopfian groups and  modules. Further sources for proofs and so on, can be found in the references of those papers. The paper by Varadarajan \cite{Var} contains many results on Hopfian and co-Hopfian objects among rings, Boolean rings, function algebras, and compact manifolds.

We did not cite Baer's paper \cite{Baer} in \cite{Me}. More glaring was our ignorance of \cite{B} and \cite{B-P}. We learned of these two papers only in the mid $2010s$. They concern {\em I-modules}, modules over a $PID$  $R$ which have submodules isomorphic to $R$. Many of the results in these papers translate easily to co-Hopfian modules. We duplicated some of them, unknowingly, in \cite{Me}.

In this article, we show how some relatively old mathematics, Matlis' concept of modules with maximal orders \cite{Eben2}, and the even older, but widely known, notion of localization of a commutative ring at a multiplicative set \cite[Ch 3]{A-M}, can inform a discussion of co-Hopfian modules. In the context of finitely generated torsion modules and direct sums of torsion cyclics, we observe similarities between the behavior of co-Hopfian modules over commutative Noetherian rings and that of the same classes of abelian groups.

Our main concern is co-Hopfian $R$-modules over a commutative Noetherian ring $R$ with $1$. However, we will present results for Hopfian and for bi-Hopfian modules when it is natural to do so.

Our first section establishes some notation and terminology. In our second section, we use Matlis' theory of modules with maximal orders to show that every submodule of $\oplus_{\mathfrak{m}\in \Omega}(\oplus_{n_{\mathfrak{m}}}E(R/\mathfrak{m}))$ is co-Hopfian given that each $n_{\mathfrak{m}}$ is a nonnegative integer, $R$ commutative Noetherian with maximal spectrum $\Omega$. In the third section, we make the elementary observation that if an $R$-module $M$ is co-Hopfian then it is necessarily an $R_S$-module, $S$ the set of $M$-regular elements of $R$. So, a co-Hopfian $R$-module is a module of fractions. If, in addition, $M$ happens to be finitely generated and faithful, then $R$ must be semilocal, and its maximal ideals must be associated primes of $M$. We can then characterize when a ring of fractions of $R$ is a {\em quoring} (equal to its own total quotient ring). We end the section with some basic results on the co-Hopficity of $P$-primary modules, $P$ a maximal ideal.  A  brief fourth section examines an instance in which we can easily determine the co-Hopfian submodules of the total quotient ring. It concludes with characterizations of finitely generated co-Hopfian modules which have finite length. Section 5 detemines the co-Hopfian submodules of the indecomposable injective module $E=E(R/P)$, $P$ an arbitrary prime of a commutative Noetherian ring $R$, and then does the same for an arbitrary co-Hopfian injective. Our sixth section  determines the Hopfian submodules of an arbitray co-Hopfian injective. In Section 7, we determine the commutative Noetherian rings for which $E(R/P)$ has ``few'' co-Hopfian submodules. This section is mainly an exercise in commutative algebra.

\section{Preliminaries}

Unless stated otherwise, $R$ is a commutative Noetherian ring with $1 \neq 0$. Subrings have the same unit as the parent ring, and ring homorphisms map identity to identity.  All $R$-modules are unital. All groups are abelian unless specified otherwise. We use $E$ to denote an arbitrary  injective $R$-module, and $E(M)$ to denote the {\em injective envelope} ({\em hull}, if you prefer) of an $R$-module $M$. As in \cite{fuchs},  $\mathbb{Z}(n)$ represents the cyclic group of order $n,$ and $\mathbb{Z}(p^{\infty})$  the Pr\"{u}fer $p$-group. We say an $R$-module $M$ is {\em $P$-primary} if $\mathrm{Ass}(M)=\{P\}$  (some prefer  {\em $P$-co-primary} \cite{Eisenbud}, \cite{Mats}). We use our terminology to be consistent with Matlis' terminology in \cite{Eben2} and the terminology of abelian groups. Every element of a $P$-primary module $M$ is annihilated by some power of $P$. If $M$ is finitely generated, then $M$ itself  is annihilated by a power of $P$ \cite[Prop 3.9]{Eisenbud}.

A {\em regular element} is an element of $R$ that is not a zerodivisor. The {\em total quotient ring} of $R$ is the ring $Q=S^{-1}R$ with $S$ the set of regular elements of $R$. For convenience, if $S$ is any multiplicative subset of $R$, we will use the notation $R_S$ as in \cite{Kap},\cite{Mats},\cite {Z-S}, and the usual $R_P$ if $S=R-P$, $P$ a prime ideal of $R$ (similarly for $R$-modules $M$).  If every regular element of $R$ is a unit, then $R=Q$. We call such a ring a {\em quoring} \cite[\S 10.7]{Faith} (the term {\em classical} is also used \cite{Lam}).  For example, Artinian rings are quorings.

If $M$ is an $R$-module, $\mu_r$ is the $R$-module endomorphism of $M$ that is left multiplication by the element $r \in R$. Note that $M$ may be the $R$-module $R$. It is well-known that $\mu_r$ is an injection if and only if $r$ is not a zerodivisor on $M$ ($r \in R, x \in M$ and $rx=0$ imply $x=0$), and a surjection if and only if $rM=M$. The $r$ which are not zerodivisors on $M$ are called {\em $M$-regular} and are those elements of $R$ that do not belong to any associated prime of $M$ \cite[p 38]{Mats}.

It is not hard to see that if $M=R$, then $\mu_r$ is a surjection if and only if $r$ is a unit (for then $\mu_r(x)=rx=1$ for some $x \in R$ \cite[Prop 8]{Hiremath}). So,  $R$ is {\em always} a Hopfian $R$-module, and is co-Hopfian if and only if it is a quoring (\cite[Prop 1.4, Remark 1.6(b)]{Var} and, independently, \cite{Me}, remark after Theorem 1.1;  the expository note \cite{Me2} provides a noncommutative generalization). It is straightforward to show that if $n$ is a positive integer, then $R^n$ is a Hopfian $R$-module \cite[Thm 14, Cor 15]{Hiremath}, and co-Hofian if and only if $R$ is a quoring (\cite[Thm 3.2]{Var} and, by a different argument \cite[Thm 2.5]{Me}).

There are some basic facts about Hopficity/co-Hopficity which we will use, often without mention. Any module satisfying the ascending chain condition ($ACC$) is Hopfian and any module satisfying the descending chain condition ($DCC$) is co-Hopfian (see \cite[Exercises 1.7, 1.10]{Lam2} or \cite[Prop 8, p 22]{Lambek}); recall that a module with $ACC$ is Noetherian, one with $DCC$ Artinian; these facts also appear as exercises in various texts, for example \cite[Exercise 1, Ch 6]{A-M}. A module of finite length has both $ACC$ and $DCC$ and so is bi-Hopfian, for example, any finitely generated module over an Artinian ring.

For any cardinal $n$, finite or infinite, we use $\oplus_nM$ to denote the direct sum of $n$ copies of the module $M$. We often use ``='' (equality) synonymously with ``$\cong$'' (isomorphism).

For convenience, we include several results (Theorems A-F) often referred to in the main text. We list Theorems A-C here, along with an illustrative example. Theorems D-F appear in the next section.

\begin{tA}[\cite{Me} Thms 1.2-1.4]
Let $R$ be a ring.
\begin{enumerate}
\item Let $M$ be a co-Hopfian (Hopfian) $R$-module. If $M$ decomposes as a direct sum of a family $\{M_i\}$ of nontrivial $R$-modules, then each $M_i$ is co-Hopfian (Hopfian).
\item Let $\{M_i\}$ be a family of nontrivial $R$-modules for which $\mathrm{Hom}(M_i,M_j) = \mathrm{Hom}(M_j,M_i)=0$ whenever $i \neq j.$ If each $M_i$ is co-Hopfian (Hopfian), then so is $\oplus M_i.$
\item   If $A,B$ are (co)-Hopfian  and either $\mathrm{Hom}(A,B)=0$ or $\mathrm{Hom}(B,A)=0,$ then $A \oplus B$ is (co)-Hopfian.
\item  Let $M$ be a co-Hopfian (Hopfian) $R$-module. If $M$ decomposes as a direct sum of a family $\{M_i\}$ of nontrivial $R$-modules, then there are only finitely many summands isomorphic to a given $R$-module $N.$
\end{enumerate}
Item (3) is due to Goldsmith and Gong \cite{G-G3}. It can be useful in dealing with direct sums, and is a special case of a more general result they established in \cite{G-G}, a result which extends easily to modules. \qed 
\end{tA}

Certain changes of ring preserve the co-Hopfian (Hopfian) property.

\begin{tB}[\cite{Me3} Thm 1.3]
Let $R \rightarrow R^{\prime}$ be an epimorphism in the category of  rings. Let $M,N$ be $R^{\prime}$-modules. Then $\mathrm{Hom}_{R^{\prime}}(M,N)=\mathrm{Hom}_R(M,N).$ \qed

\end{tB}

This result is contained in Proposition XI.1.2 in \cite{Bo}. It can also be found  at the Stacks Project (stacks.math.columbia.edu/tag/04VM) in the form as written. We will use this theorem in the particular cases in which $R^{\prime}$ is a quotient of $R,$ a localization of $R,$ or, if $R$ is a domain, the fraction field of $R.$ The modules $M$ and $N$ are $R$-modules by restriction of scalars.

Note that in the particular case in which $M=N$, Theorem B implies that $M$ is co-Hopfian (Hopfian) as an $R$-module if and only if it is co-Hopfian (Hopfian) as an $R'$-module.

\begin{tC}[\cite{Me3} Thm 2.1]
Let $E \cong \bigoplus_{P\in \mathrm{Ass}(R)}(\oplus_{n_P}E(R/P))$ be an injective $R$-module. Then $E$ is co-Hopfian if and only if $n_P$ is finite for all $P.$
\end{tC}

In \cite[Cor 2.5]{Me3} we noted that if we restrict the sum in the theorem to minimal primes $P,$ then $E$ is Hopfian, so bi-Hopfian. If $R$ is reduced Noetherian, then $E$ is bi-Hopfian if and only if the sum is restricted to minimal primes \cite[Cor 2.6]{Me3}.

Theorem A suggests that if we have a satisfactory structure theory for a class of modules, then we can try to exploit that structure to determine which modules in the class are Hopfian or co-Hopfian. 

Let $G$ be an abelian group whose torsion subgroup $tG$ is cotorsion. Then $G=tG\oplus G/tG$ \cite[\S 54]{fuchs}. Hence, $G$ is co-Hopfian if and only if both $tG$ and $G/tG$ are (Theorem A, or, see the remark following the proof of \cite[Thm 3.1]{Me}). 

The torsion-free group $G/tG$ is co-Hopfian if and only if $G/tG \cong \mathbb{Q}^n$ for some nonnegative integer $n$ (torsion-free + co-Hopfian = divisible \cite[Thm 2.1]{Me}). A torsion cotorsion group is the direct sum of a divisible group $D$ and a bounded group $B$  \cite[Cor 54.4]{fuchs}, that is $nB=0$ for some positive integer $n$, which we will always assume is the least such. A bounded group $B$ is the direct sum of torsion cyclic groups \cite[Thm 17.2]{fuchs}. But $B$ is the direct sum of finitely many $p$-components $B_p$, one for each prime divisor of $n$, and each of these is a direct sum of cyclics. Now $B_p$ is reduced \cite[p 98]{fuchs}, and at most countable (Theorem A). So, $B_p$ is a finite group \cite[Prop 2.2(vi)]{G-G}. Finally, the divisible group $D=\oplus_p(\oplus_{n_p}\mathbb{Z}(p^{\infty}))$, the sum over the positive primes (Theorem C; divisible = injective for abelian groups \cite[Thm 24.5]{fuchs}). Note that no $\mathbb{Z}(p^{\infty})$ is Hopfian. 

\begin{thm}
Let $G$ be an abelian group for which $tG$ is cotorsion. Then $G$ is co-Hopfian if and only if  $G \cong \mathbb{Q}^n \oplus [\oplus_p(\oplus_{n_p}\mathbb{Z}(p^{\infty}))] \oplus B$, with $n$ and the $n_p$ nonegative integers and $B$ a finite group. Moreover, $G$ is Hopfian if and only if $G \cong \mathbb{Q}^n \oplus B$, $n$ a positive integer and $B$, again, a finite group.\qed
\end{thm}

\section{Maximal orders; Localization}

Matlis \cite{Eben2} introduced modules with maximal orders as  natural generalization of modules with $DCC$. An $R$-module $M$ has {\em maximal orders} if for each nonzero $x \in M$, the only prime ideals containing $\mathrm{ann}_R(x)$ are maximal ideals. For example, torsion modules over a $1$-dimensional Noetherian domain, specifically torsion abelian groups, are such modules.

An $R$-module $M$ with maximal orders has, for each $\mathfrak{m}\in \Omega$, the maximal spectrum of $R$, an $\mathfrak{m}$-primary component $X_{\mathfrak{m}}(M)$  defined by
$$X_{\mathfrak{m}}(M) = \{x \in M : \mathfrak{m}^n x = 0,\hspace{.3em} \mbox{for some positive integer $n$}\}.$$
It is clear that $\mathfrak{m}$ is the only associated prime of $X_{\mathfrak{m}}(M)$.

The following theorems give the portions of three of Matlis' theorems  most relevant to us. There are obvious similarities between modules with maximal orders and torsion abelian groups. In particular, the third theorem shows that there is a ``primary decomposition'' for modules with maximal orders. We have changed the notation slightly to conform with our current notation.

\begin{tD}[\cite{Eben2}, Prop 2]
Let $\mathfrak{m}$ be a maximal ideal of $R$, and let $\hat{R}_{\mathfrak{m}}$ be the $\mathfrak{m}$-adic completion of $R_{\mathfrak{m}}.$ Then if $M$ is an $\mathfrak{m}$-primary $R$-module, we have natural isomorphisms:
$$M \cong M \otimes_R R_{\mathfrak{m}} \cong M \otimes_R \hat{R}_{\mathfrak{m}}$$.
\end{tD}

\begin{tE}[\cite{Eben2}, Prop 3]
Let $M$ be an $R$-module. Then the following are equivalent:
\begin{enumerate}
\item $M$ has DCC.
\item $M$ is a submodule of $E_1 \oplus \ldots \oplus E_n,$ where $E_i = E(R/\mathfrak{m}_i),$ $\mathfrak{m}_i$ a maximal ideal of $R.$
\item $M$ has maximal orders and finitely generated socle. \qed
\end{enumerate}
\end{tE}

The {\em socle} of $M$ is the submodule $\mathrm{soc}(M)$ of $M$ that is the sum of the simple submodules of $M$. In fact, $\mathrm{soc}(M)$ is a semisimple module \cite[Prop 3.2]{S-V}.

\begin{tF}[\cite{Eben2}, Thm 1]
Let $M$ be an $R$-module. Then the following are equivalent:
\begin{enumerate}
\item $M$ has maximal orders.
\item $M$ is an essential extension of its socle.
\item $M$ is a submodule of $\bigoplus _{\mathfrak{m} \in \Omega}(\oplus_{n_{\mathfrak{m}}}E(R/\mathfrak{m}))$, for cardinal numbers $n_{\mathfrak{m}}.$
\item $M =  \bigoplus_{\mathfrak{m} \in \Omega} X_{\mathfrak{m}}(M).$
\item Every finitely generated submodule of $M$ has finite length. \qed \\
\end{enumerate}
\end{tF}

The following corollary to item 4 of Theorem F reduces questions of co-Hopficity or Hopficity for an $R$-module $M$ having maximal orders to the corresponding questions for the primary components of $M$.

\begin{cor}
Let $M$ be an $R$-module having maximal orders. Then $M$ is co-Hopfian (resp. Hopfian) if and only if $ X_{\mathfrak{m}}(M)$ is co-Hopfian (resp. Hopfian) for each $\mathfrak{m} \in \Omega$.
\end{cor}

\begin{proof}
If $M$ is co-Hopfian (resp. Hopfian), then so is each $ X_{\mathfrak{m}}(M)$ (Theorem A). Conversely, if each $ X_{\mathfrak{m}}(M)$ is co-Hopfian (resp. Hopfian), then so is $M$  since $\mathrm{Hom}(X_{\mathfrak{m}}(M),X_{\mathfrak{n}}(M))=0$ if $\mathfrak{m} \neq \mathfrak{n}$ (Theorem A again).
\end{proof}

We can specify immediately which submodules of a co-Hopfian injective $R$-module $E$ are co-Hopfian, given that $E$ has maximal orders.

\begin{thm}
Let $E = \bigoplus_{\mathfrak{m}\in \Omega}(\oplus_ {n_{\mathfrak{m}}}E(R/\mathfrak{m})),$ with each $n_{\mathfrak{m}} < \infty.$ Any submodule $M$ of $E$ is co-Hopfian. 
\end{thm}

\begin{proof}
We know $E$ is co-Hopfian (Theorem C). Following Corollary 2.1, let $X_{\mathfrak{m}}(M)$ be a primary component of $M$. By definition, $X_{\mathfrak{m}}(M) \subseteq \oplus_{n_{\mathfrak{m}}} E(R/\mathfrak{m})$. Since $n_{\mathfrak{m}}$ is finite, $X_{\mathfrak{m}}(M)$ has DCC (Theorem E) and so is co-Hopfian. Hence,  by Theorem A, $M=\oplus X_{\mathfrak{m}}(M)$ is co-Hopfian, per the proof of Corollary 2.1.
\end{proof}

\begin{thm}
Let $M$ be an $R$-module having maximal orders.
\begin{enumerate}
\item If $E(M)$ is co-Hopfian, then so is $M$.
\item If $M$ is finitely generated, then $M$ is bi-Hopfian.
\item If $M$ is finitely generated, then every submodule of $M$ is bi-Hopfian.
\item  If $M =  \bigoplus_{\mathfrak{m} \in \Omega} X_{\mathfrak{m}}(M)$ and $X_{\mathfrak{m}}(M)$ is finitely generated for each $\mathfrak{m} \in \Omega$, then $M$ is bi-Hopfian.
\item If $M =  \bigoplus_{\mathfrak{m} \in \Omega} X_{\mathfrak{m}}(M)$, then $X_{\mathfrak{m}}(M)$ is the localization of $M$ at $\mathfrak{m}$.
\end{enumerate}
\end{thm}

\begin{proof}
(1): Immediate, since each associated prime of $M$ is the annihilator of a nonzero $x \in M$, hence a maximal ideal. So, if $E(M)$ is co-Hopfian, then it has the form of the module $E$ of Theorem 2.2.

(2), (3): Both follow directly from Theorem F.

(4): Follows from (3) and Corollary 2.1.

(5): The localization of $M$ at $\mathfrak{m}$, $M_{\mathfrak{m}}$, is (isomorphic to) $R_{\mathfrak{m}}\otimes M$. Tensor product distributes over direct sum \cite[Thm 8.87]{Joe2}, so we need only know $R_{\mathfrak{m}}\otimes X_{\mathfrak{n}}(M)$ for maximal ideals $\mathfrak{m},\mathfrak{n}$. If $\mathfrak{m}=\mathfrak{n}$ this product is $X_{\mathfrak{n}}(M)$ (Theorem C). If $\mathfrak{m}\neq \mathfrak{n}$, the product is zero. If $x$ is a nonzero element of $X_{\mathfrak{n}}(M)$, then $\mathfrak{n}^kx=0$ for some positive integer $k$. Let $a\in \mathfrak{n}\setminus \mathfrak{m}$. Then $a^k \in \mathfrak{n}^k$, and $a^kx=0$. So,
$$\frac{r}{s} \otimes x= \frac{ra^k}{sa^k} \otimes x = \frac{r}{sa^k} \otimes a^kx=\frac{r}{sa^k} \otimes 0 = 0.$$
\end{proof}

In Section 3, we will have more to say about Theorem 2.3. Theorem 2.2 generalizes the co-Hopfian part of an observation made in \cite[Rem 3.6]{Me3}. There we showed that if $R$ is Artinian, then the module $E$ of Theorem 2.2 and all of its submodules are bi-Hopfian. If $R$ is not Artinian, then $E(R/\mathfrak{m})$ may not be Hopfian. For example, $\mathbb{Z}(p^{\infty})=E(\mathbb{Z}/(p))$ is not Hopfian (see the discussion in \cite[\S 2]{Me3}). \\

The localization of a ring $R$ at a prime ideal $P$ is used widely in commutative algeba and in algebraic geometry. More generally, we have the notion of localization of $R$ at a multiplicatively closed subset $S$ of $R$. The construction extends directly to $R$-modules (see, for example,  \cite[Ch 3]{A-M}). For an appropriate choice of $S$, localizing an $R$-module $M$ at $S$ is useful in studying co-Hopfian modules, and enters our discussion in the guise of scalar multiplication. Recall that scalar multiplication is definitive in determining if a ring $R$ is co-Hopfian as an $R$-module (Section 1).

 For general $R$-modules $M$ the situation is more complex. The $R$-module endomorphism $\mu_r$ of $M$ is an injection if and only if $r$ is $M$-regular, but $\mu_r$ can be a surjection without $r$ being a unit, as with multiplication by $2$ on the abelian group $\mathbb{Q}$ (this situation cannot occur if $M$ is finitely generated and faithful, see Lemma 2.5). We will see that, at least as far as finitely generated modules are concerned, the interplay between certain primes in $\mathrm{Ass}_R(M)$, the set of associated primes of the $R$-module $M$,  and the ideals of $R$ is important.

We begin with an elementary observation. The set $S$ of $M$-regular elements of $R$ is multiplicatively closed. If $M$ is co-Hopfian, then $\mu_r$ must be an automorphism of $M$ for each $r$ in $S$, precisely the condition for $M$ to be an $R_S$-module \cite[\S 6.5]{Reid}, \cite[Prop II.3.7]{Bo}. So, our first result is clear.

\begin{thm} Let $M$ be an $R$-module and $S$ the set of $M$-regular elements of $R.$ If $M$ is co-Hopfian, then $M$ is an $R_S$-module and $M \cong M_S$ .\qed
\end{thm}

So, a co-Hopfian $R$-module $M$ is a module of fractions over a ring of fractions. Matlis \cite[Lemma 3.2]{Eben} used this approach, at least implicitly, to show that the indecomposable injective module $E(R/P)$ is an $R_P$-module, $P$ a prime ideal of a Noetherian ring $R$.

We can say more if $M$ is a finitely generated faithful $R$-module. It is well-known that if $I=\mathrm{ann}_R(M),$  then $M$ is a faithful $R/I$-module \cite[p 20]{A-M}. Since $M$ is co-Hopfian as an $R$-module if and only if it is co-Hopfian as an $R/I$-module (Theorem B), there is no loss of generality in assuming that $M$ is a faithful $R$-module.  

Suppose, in addition, that $M$ is finitely generated. Let $U$ be the set of units of $R.$ If $r \in U,$ it is clear that $\mu_r$ is an automorphism of $M$. Conversely, if $\mu_r$ is an automorphism of $M$, then $JM=M,$ where  $J=(r).$  Hence, there is $x \in R$, $x \equiv 1 (\mathrm{mod} \: J),$ for which $xM=0$ \cite[Cor 2.5]{A-M}. Since $M$ is faithful, $1 \in J$ (since $x$ must be $0$), so $J=R$. Thus, $r \in U$. We summarize in a lemma.

\begin{axe}
If $M$ is a finitely generated faithful $R$-module and $r \in R$, then $\mu_r$ is an automorphism of $M$ if and only if $r$ is a unit.\qed
\end{axe}

The lemma provides sufficient leverage for the following result. 

\begin{thm}
Let $R$ be Noetherian and $M$ a finitely generated faithful $R$-module. If $M$ is co-Hopfian, then the maximal elements of $\mathrm{Ass}_R(M)$ are the only maximal ideals of $R$, so that $R$ is semilocal. In addition, the $M$-regular elements of $R$ are regular elements of $R$.
\end{thm}

\begin{proof}
Let $P_1, \ldots , P_k$ be the distinct maximal elements of the finite set $\mathrm{Ass}_R(M)$. Renumbering, if necessary, suppose $P_1$ is not a maximal ideal. Then $P_1 \subsetneq \mathfrak{m}$ for some maximal ideal $\mathfrak{m}$. Since the $P_i$ have no containments between them, the ideal $\mathfrak{m}$ is not equal to any of the $P_i$, including $P_1$. In particular, $\mathfrak{m}$ is not contained in any of the $P_i$. Prime avoidance gives us an $r \in \mathfrak{m}$ that is contained in no $P_i$. Hence, $r$ is $M$-regular but not a unit of $R$. Thus, $\mu_r$ is an injective endomorphism of $M$ which is not an isomorphism, so $M$ is not co-Hopfian.

A similar argument shows that the $P_i$ are the only maximal ideals of $R$. If $R$ were to have a maximal ideal $\mathfrak{m}$ other than the $P_i$, prime avoidance once again gives an $M$-regular element that is not a unit of $R$.

The last statement follows since the $M$-regular elements are units in $R$. 
\end{proof}

\begin{remk}
If $M$ has nonzero annihilator $I$, the theorem applies to $A=R/I$ and implies that $I$ is contained in only finitely many maximal ideals of $R$ (just as in the case of finitely generated torsion groups).

The theorem implies that if a Noetherian $R$ admits a finitely generated faithful co-Hopfian modules, then $R$ must be semilocal.

To us, the most interesting part of Theorem 2.6 may be the last statement, which is reminiscent of Auslander's Zero Divisor Conjecture (now a theorem, see \cite{Naseh} for details). It is independent of Auslander's condition that $M$ have finite projective dimension. \qed \\
\end{remk}

\noindent{\bfseries Further remarks.} (1):  Fuchs \cite[p 20, Examples 2,3]{fuchs} discusses Theorem D in the context specific to $p$-groups. Kaplansky \cite[p 42]{KapAG} notes the special case for modules over a $DVR$. 

(2): Theorem 2.4 is useful in showing a module is {\em not} co-Hopfian. For example, let $R$ be an integral domain with fraction field $Q$ and $M$ a torsion-free $R$-module. If $M$ is not a $Q$-vector space, then $M$ is not co-Hopfian (compare with \cite[Thm 3]{B-P}). There is a partial converse: if $M$ is a finite-dimensional $Q$-vector space, then $M$ is co-Hopfian. We showed this in \cite[Thm 2.4]{Me} using a different argument than \cite{B-P}.

Note further that if $N$ is an $R$-submodule of the $R$-module $M$ and is co-Hopfian, then $N$ must also be an $R_S$-module ($M$-regular elements are $N$-regular since $\mathrm{Ass}(N) \subseteq \mathrm{Ass}(M)$), an observation we exploit in Sections 4 and 5. 

(3): The converse of Theorem 2.4 is not true generally. As Fuchs pointed out, a $p$-group is naturally a $\mathbb{Z}_{(p)}$-module, but some are co-Hopfian (any finitely generated one), while others are not (for example, a countable direct sum of copies of $\mathbb{Z}(p)$ (Theorem A)). 

The converse can fail also if $R_S$ is isomorphic to $R$. For example,  let  $G=\prod_p \mathbb{Z}(p^{\infty}),$ the product being over the positive primes. Then $G$ has torsion subgroup $tG = \oplus_p  \mathbb{Z}(p^{\infty})$, which is a direct summand of $G$ (it is divisible).  It is also co-Hopfian (Theorem A). The complementary summand is isomorphic to a direct sum of uncountably many copies of $\mathbb{Q}$ \cite[Exercise 9.36]{Joe}, and so is not co-Hopfian (Theorem A). Therefore, $G$ is not co-Hopfian (again, Theorem A). Since every nonzero prime ideal of $\mathbb{Z}$ is an associated prime of both $G$ and $tG,$ the sets $S$ and $S^{\prime}$ of $G$-regular and $tG$-regular elements, respectively, consist of just $\pm 1.$ Hence, $\mathbb{Z}_S$ and $\mathbb{Z}_{S^{\prime}}$ are both isomorphic to $\mathbb{Z}.$ But, $tG$ is co-Hopfian and $G$ is not.

The late H. Storrer was kind enough to remind us of this example when he pointed out that the direct product portion of \cite[Thm 1.3]{Me} is not correct.

(4): The converse of the Theorem 2.6  is not true. Let $(R, \mathfrak{m})$ be a local domain of positive Krull dimension.  Let $P$ be a  prime ideal properly contained in $\mathfrak{m}.$ Let $M = R/\mathfrak{m} \oplus R/P.$ Then $\mathrm{Ass}(M)=\{\mathfrak{m},P\}$ and $\mathrm{ann}(M)=P.$ So, $M$ is a nonzero finitely generated faithful $R/P$-module, $R/P$ is local, and $\mathfrak{m}/P$ is the maximal element in $\mathrm{Ass}_{R/P}(M).$ However, $M$ is {\em not} a co-Hopfian $R$-module since $R/P$ is not co-Hopfian: any element of $\mathfrak{m} \setminus P$ is regular on $R/P$ but not a unit of $R.$ The module $M$ is an example of a finitely generated $R$-module for which $\mathrm{Ass}_R(R/\mathrm{ann}(M))$ is properly contained in $\mathrm{Ass}_R(M)$ (\cite[Exercise 3.40D]{Lam}, \cite[p 85]{Lam2}). Vasconcelos \cite{Wolmer} used a similar setup to show that if a commutative ring $R$ has positive Krull dimension, then $R$ has a finitely generated module, the cyclic module $R/P$, $P$ a nonmaximal prime, that is not co-Hopfian. 

(5): Any Noetherian ring $R$, semilocal or not, supports a finitely generated faithful $R$-module, namely $R$ itself. This module is co-Hopfian if and only if $R$ is a quoring. For $n$ a positive integer, $R^n$ is also finitely generated, faithful, and co-Hopfian. We have many semilocal quorings, the total quotient rings of Noetherian rings.

No  integral domain $R$ supports a finitely generated faithful torsion module. If $x_1,x_2, \ldots ,x_n$ are generators for such a module $M$, then for each $i$ there is a nonzero $a_i$ in $R$ for which $a_ix_i=0$. But then $a=\prod a_i\neq 0$ and $aM=0$.

If $R$ is semilocal, but not a field, let us say that an element $x$ of an $R$-module $M$ is a torsion element if $rx=0$ for some nonzero $r \in R$ (this is {\em not} the standard definition). It makes sense to call an $R$-module {\em mixed} if it has both torsion and non-torsion elements. For example a ring $R$ having zero divisors is mixed since units are not torsion.

 Let $R$ be Artinian with distinct maximal ideals $\mathfrak{m}_1, \ldots , \mathfrak{m}_k.$ In this case, the Jacobson radical of $R$ and the nilradical, $\mathrm{Nil}(R)$, of $R$ coincide and are equal to $\cap_{i=1}^k \mathfrak{m_i}$ \cite[Cor 8.2]{A-M}. Moreover, $\mathrm{Nil}(R)$ is nilpotent. Since the $\mathfrak{m_i}$ are maximal (and distinct), $\cap \mathfrak{m_i}=\prod \mathfrak{m_i}$. The nilpotence of $\mathrm{Nil}(R)$ implies that $\prod \mathfrak{m}_i^t=0$ for some positive integer $t$.
 
 Let $M=\oplus_{i=1}^k R/\mathfrak{m}_i^t.$ Note that $\mathrm{Ann}(M)=\cap \mathfrak{m}_i=\prod \mathfrak{m}_i=0$, so that $M$ is a finitely generated and faithful $R$-module. Also, $M$ is co-Hopfian (Theorem A(2)) and torsion. It is immediate that $M=\oplus R/\mathfrak{m}_i^{n_i}$ has the same properties provided that $n_i \geq t$ for all $i$. The mixed $R$-module  $$M=R^n \oplus (\oplus_{i=1}^k R/\mathfrak{m}_i^{n_i})$$
 is finitely generated and faithful. It is co-Hopfian for any positive integer $n$ provide that each $n_i \geq t$ ($M$ is finitely generated). \qed

 \section{Finitely generated modules}

 We recall that a finitely generated $R$-module $M$ is Hopfian. So, if such a module is co-Hopfian, it is necessarily bi-Hopfian. The simplest finitely generated modules are the cyclic modules. Every cyclic $R$-module is of the form $R/I$, for an ideal $I$ of $R$ \cite[Prop 7.12]{Joe2}. But $R/I$ is co-Hopfian as an $R$-module if and only if it is co-Hopfian as an $R/I$-module (Theorem B) if and only if the {\em ring} $R/I$ is a quoring. Note that the ring $R/I$ is the endomorphism ring of the module $R/I$. Our first result determines when the ring $R/I$ is a quoring. It is actually a corollary to Theorem 2.6.

\begin{thm}
Let $R$ be Noetherian and $I$ an ideal of $R$.  Then $R/I$  is a quoring if and only if it is semilocal and each of its maximal ideals is an associated prime of the $R/I$-module $R/I$. The case $I=0$ characterizes Noetherian quorings.
\end{thm}

\begin{proof}
($\Rightarrow$) A quoring is co-Hopfian so Theorem 2.6 verifies this implication.

($\Leftarrow$) The condition on maximal ideals implies each element of $R$ is a unit.
\end{proof}

Within the cyclic modules, the easiest to analyze are the $R/I$ for $I$ a $P$-primary ideal, $P$ a prime of $R$. But, $R/I$  is $P$-primary if and only if $I$ is $P$-primary \cite[Thm 7.8]{Reid}. Before investigating these modules, we pause for a result about finitely generated modules that will be useful in our analysis.

Let $M$ be a finitely generated $R$-module with annihilator $I$. Again, we see that if $M$ is $P$-primary, then $I$ is $P$-primary. Generally speaking, $\mathrm{Ass}(R/I) \subseteq \mathrm{Ass}(M)$. If $R/I \neq 0$, then $\mathrm{Ass}(R/I) \neq \emptyset$. But, $\mathrm{Ass}(M)=\{P\}$, so $I$ is $P$-primary (the converse is not true; see the ``Further remarks'' to Section 2). If $M$ is finitely generated and $P$-primary, then $M$ is co-Hopfian if and only if $P$ is maximal.

\begin{thm}
Let $P$ be a prime ideal of $R$ and $M$ a finitely generated $P$-primary $R$-module with annihilator $I$. The following are equivalent.
\begin{enumerate}
\item $P$ is a maximal ideal of $R$.
\item $R/I$ is an Artinian ring.
\item  $M$ is a co-Hopfian $R$-module.
\end{enumerate}
Consequently, if $P$ is maximal, then every nonzero submodule of $M$ is co-Hopfian, but if $P$ is not maximal, no nonzero submodule of $M$ is co-Hopfian.
\end{thm}

\begin{proof}
 Since $\mathrm{Ass}(M)=\{P\}$, $P$ is the minimal prime containing $I$ \cite[Thm 3.1]{Eisenbud}. So, $P = \mathrm{rad}(I)$  by definition of radical of an ideal. 

(1) $\Rightarrow$ (2): Since $R$ is Noetherian, $P^n \subseteq I$ for some positive integer $n$ \cite[Prop 7.14]{A-M}. Hence, $P^n$ annihilates $M.$ But $P$ maximal implies $R/I$ is Artinian \cite[Cor 2.17]{Eisenbud}.

(2) $\Rightarrow$ (3):  $M$ is a finitely generated $R/I$-module, so a co-Hopfian $R/I$-module, hence a co-Hopfian $R$-module (Theorem B). 

(3) $\Rightarrow$ (1): $\mathrm{Ass}(M)=\{P\}$. If $P$ is not maximal, it is contained in some maximal ideal $\mathfrak{m}$. Any $r$ in $R \setminus \mathfrak{m}$ yields an endomorphism $\mu_r$ of $M$ that is an injection but not a surjection, whence $M$ is not co-Hopfian.

To conclude, note that every nonzero submodule of $M$ is finitely generated and $P$-primary, so the equivalence of (1) and (3) yields the final statements.
\end{proof}

\begin{cor}
\begin{enumerate}
\item Let $I$ be a $P$-primary ideal of $R$. Then $(R/I)^n$ is co-Hopfian if and only if $P$ is maximal and $n$ is a positive integer.
\item Let $\{I_j\}$, $1 \leq j \leq n,$ be $P$-primary ideals. Then $M=\oplus_{j=1}^nR/I_j $ is co-Hopfian (in fact, bi-Hopfian) if and only if $P$ is maximal.
\item Let $J$ an index set, $\{M_j\}$ $R$-modules with $M_j$ finitely generated and $P_j$-primary. Suppose further that $P_i \neq P_j$ if $i \neq j.$ Then $M=\oplus M_j$ is co-Hopfian (in fact, bi-Hopfian) if and only if each $P_j$ is maximal.
 \item Let $I$ be a $P$-primary ideal of $R$. If $P$ is maximal, then every submodule of $R/I$ is co-Hopfian. If $P$ is not maximal, then no nonzero submodule of $R/I$ is co-Hopfian. 
\end{enumerate}
\end{cor}

\begin{proof}
(1): Immediate from the theorem since $(R/I)^n$ is finitely generated and $P$-primary. 

(2): $M$ is finitely generated and $P$-primary, so $M$ is co-Hopfian if and only if $P$ is maximal. But $M$ is also Hopfian since it is finitely generated over a commutative ring (Section 1).

(3): If $M$ is co-Hopfian, then each $M_j$ is (Theorem A). Hence, each $P_j$ is maximal.  If each $P_j$ is maximal, then $\mathrm{Hom}(M_i,M_j)=0$ if $i \neq j$ (see Lemma 5.14). So $M$ is co-Hopfian (Theorem A). As in (2), each $M_j$ is Hopfian. Thus, each $M_j$ is bi-Hopfian, so $M$ is bi-Hopfian (Theorem A).

(4): $\mathrm{Ass}(R/I) = \{P\}$ (\cite[Thm 7.8]{Reid}, \cite[Thm 6.6]{Mats}), so $R/I$ is $P$-primary and finitely generated. Thus, the final statements of the theorem apply.
\end{proof}

\begin{remk}
A particular case of item (2) in the corollary is worth noting. If $I=P$ is a maximal ideal, then the positive integral powers $P^n$ of $P$ are $P$-primary, so any $R/P^n$ is co-Hopfian (bi-Hopfian, even), as is any finite direct sum of such modules for various $n.$ This is a situation reminiscent of finitely generated torsion abelian $p$-groups.

Specifically, if $\Omega$ is the maximal spectrum of $R$, then $M =\oplus_{P \in \Omega}(\oplus_{j=1}^t R/P^{n_j})$ is a co-Hopfian $R$-module, where for each $P$ the $n_j$ are various positive integers.  This follows directly from item (3) of the corollary. Observe that a direct sum of cyclics $R/P^n$ for arbitrary primes $P$ cannot be co-Hopfian if any $P$ is not maximal (Corollary 3.3).  

Note that we can replace $\oplus_{j=1}^t R/P^{n_j}$ by $\oplus_{j=1}^t R/Q_j$ for various $P$-primary ideals $Q_j.$ The annihilator of $\oplus_{j=1}^t R/Q_j$ is the $P$-primary ideal $\cap Q_i,$ so $\oplus_{j=1}^t R/Q_j$ remains $P$-primary. In fact, the theorem allows us to replace $\oplus_{j=1}^t R/Q_j$ by any finitely generated $P$-primary module, thus yielding an alternative proof of item (3) of the corollary.

Let $M=\oplus_{j=1}^t  M_j$ with $M_j$ a finitely generated $P_j$-primary module, $P_j$ a prime ideal of $R$. As above, if any $P_j$ is not a maximal ideal, then $M$ is not co-Hopfian. \qed
\end{remk}

The situation for general ideals $I$ is more delicate. The ideal $I$ has a primary decomposition  $I = \cap_{i=1}^k Q_i,$ with each $Q_i$ being $P_i$-primary for prime ideals $P_i$ of $R$. Since $R/I$ must be semilocal with each maximal ideal an associated prime of $R/I$, some of the $P_i$ must be maximal ideals of $R$ and these must contain the other $P_i$. However, if $I$ has maximal orders, then all $P_i$ are maximal. So,  $R/I = \bigoplus_{i=1}^k R/Q_i$. In this case, $R/I$ is co-Hopfian by Corollary 3.3, and $R/I$ decomposes as a direct sum of $P$-primary cyclic modules for various maximal ideals $P$, a situation reminiscent of decomposing $\mathbb{Z}(n)$ for composite $n$. If $R$ is a $1$-dimensional domain, this observation provides a complete characterization.

\begin{thm}
Let $R$ be a $1$-dimensional Noetherian domain. Then $R/I= \bigoplus_{i=1}^k R/Q_i$, $Q_i$ a $P_i$-primary ideal for some maximal ideal $P_i$, $I$  a nonzero ideal of $R$. Hence, $R/I$ is co-Hopfian. \qed
\end{thm}

The domains in Theorem 3.5 are abundant. A quick example is the coordinate ring of an irreducible affine (plane) curve. Many others can be manufactured by taking the quotient of the polynomial ring $k[x]$ over a field $k$ by a nonmaximal prime $P$ that is contained in no primes other than maximal ideals. We used to say such primes had {\em depth} $1$ (\cite[p 240]{Z-S} or \cite[Remark, p 121]{A-M}), but current usage favors {\em coheight} $1$ (coheight = dimension of $R/P$) \cite[p 30]{Mats}. Depth now refers to a deeper concept \cite[p 429]{Eisenbud}, \cite[Thm 16.7 with $M=R$]{Mats} (which used to be called {\em grade} \cite[p 89]{Kap}).

\begin{remk}
If $I$ is a nonzero ideal of a $1$-dimensional Noetherian domain $R$, then $R/I$ has finite length \cite[Thm 90]{Kap} and so is Artinian. A ring, commutative or not, whose every proper quotient is Artinian is called {\em almost Artinian}. The only commutative examples are Artinian rings (obviously), and $1$-dimensional Noetherian domains \cite{Cohen}. The latter rings are very interesting. See \cite{Eben4} for instance.   \qed 
\end{remk}

 \section{Submodules of the total quotient ring}

According to Lam \cite[Exercise 12.10]{Lam2}, E. Davis \cite{Davis} was the first to observe that a commutative ring $R$ has a semilocal total quotient ring $Q$ if and only if the set of zerodivisors of $R$ is the union of a finite number of prime ideals. In the Noetherian case, we know that the set of zerodivisors of $R$ is the union of the associated primes $P_1, \ldots , P_k$ of $R$. These ideals help to determine which $R$-submodules of $Q$ are co-Hopfian.

In its most general sense, a {\em fractional ideal} of $R$ is an $R$-submodule $I$ of $Q$ (\cite[p 31]{Lam}, or \cite[p 257]{Eisenbud}). So, any ideal of $R$ is a fractional ideal. If a fractional ideal $I$ is a co-Hopfian $R$-module, then it must be a $Q$-module (Theorem 3.1), hence an ideal of $Q$. If $s \in R$ is regular, then $s$ (actually $s/1$) is regular in $Q$ \cite[Cor 1, p 43]{Z-S}, hence a unit of $Q$. So, $\mu_s$ is an automorphism of $Q$. It restricts to an injective endomorphism of $I$. But, $I$ co-Hopfian assures $sI=I$, and so $I$ is a $Q$-module. Thus, the problem of determining which fractional ideals of $R$ are co-Hopfian $R$-modules reduces to determining which ideals of $Q$ are co-Hopfian $Q$-modules. Since $Q$ is a co-Hopfian $R$-module (it's a quoring), it suffices to consider proper ideals of $Q$. If the associated primes of $R$ are all minimal primes, it is an easy matter to determine which ideals of $Q$ are co-Hopfian (Theorem 4.3).

We begin  with a general feature of the total quotient ring of a commutative ring.

\begin{thm}
Let $R$ be commutative, but not equal to $Q$.  Then $Q$ is a minimal co-Hopfian $R$-module extension of $R$.
\end{thm}

\begin{proof}
Let $I$ be a fractional ideal of $R$ which contains $R$. If $I$ is co-Hopfian, then it must be an ideal of $Q$. But, since $1_R \in I$, we must have $I=Q$.
\end{proof}

The simplest case is when $R \neq Q$ is a domain. Then $0$ is the only associated prime of $R$ and $Q$ is a field. Since the only ideals of $Q$ are $0$ and $Q$, no proper nonzero fractional ideal of $R$ is co-Hopfian, since none of them is an ideal of $Q$.

\begin{thm}
Let $R \neq Q$ be a domain. No proper nonzero fractional ideal of $R$ is co-Hopfian.  \qed
\end{thm}

More generally, if $R$ is Noetherian, it has an Artinian total quotient ring precisely when its associated primes (those associated to the 0-ideal) are minimal primes of $R$ (\cite[Example 2, \S XV.3]{Bo},\cite[Prop 12.22]{Lam}). Noetherian domains are such rings. So is a  Noetherian ring that is reduced (= has no nonzero nilpotents), since the 0-ideal is the intersection of its minimal primes. This class of rings includes the coordinate rings of affine algebraic sets \cite[p 33]{Eisenbud}. Of course, if $Q$ is Artinian, then each of its ideals has $DCC$, and so is co-Hopfian as a $Q$-module. But then they are co-Hopfian $R$-modules (Theorem B).

\begin{thm}
Let $R$ be a Noetherian ring whose associated primes are all minimal primes. Every ideal of its total quotient ring is co-Hopfian. \qed
\end{thm}

Let $R$ be a Noetherian domain and $M$ a torsion-free $R$-module. In \cite[Prop 11]{Hiremath}, Hiremath shows that if $E(M)$ is Hopfian, then so is $M$. Theorem 4.2 shows that such a result is not possible in the co-Hopfian case, at least for rank one torsion-free $M$. In general, if $M$ is torsion-free, so is $E=E(M)$. Hence, $E$ is a direct sum of copies of $Q$, the fraction field of $R$, the number of summands being equal to the rank of $M$. But then $E$ is co-Hopfian if and only if the number of copies is finite. Arguing as above, $M$ is co-Hopfian if and only if it is a $Q$-vector subspace of $E$. So, if $E=Q^n$, the only co-Hopfian submodules of $E$ (up to isomorphism) are $0$ and the $Q^k,$ $1 \leq k \leq n.$ However, $Q$ is the injective envelope of any nonzero proper fractional ideal $I$ of $R$, of which there are many, and none of these is a $Q$-vector space by Theorem 4.1. So, we can build many $R$-submodules of $E$ that are not co-Hopfian, even though $E$ is.

\begin{thm}
Let $R$ be a domain with quotient field $Q$. If $E$ is a torsion-free injective $R$-module, then $E$ is co-Hopfian if and only if it has finite rank, in which case its only co-Hopfian submodules are its $Q$-vector subspaces.  \qed
\end{thm}

\begin{remk}
The tidy Exercise 9.7 of \cite{Joe2} contains useful information about when an $R$-module $M$ is a $Q$-vector space, $R$ an integral domain. \qed
\end{remk}

We note that the total quotient ring is useful in answering an exercise that E. Enochs once suggested to us: if $R$ is commutative Noetherian and $M$ is a finitely generated co-Hopfian $R$-module, must $M$ have finite length? If $R$ is Noetherian, so is its total quotient ring $Q$, and $Q$ is a finitely generated co-Hopfian $Q$-module. However, $Q$ is an Artinian ring, and so Artinian as a $Q$-module, if and only if the associated primes of $R$ are exactly the minimal primes of $R$ (\cite[Prop 12.22]{Lam} or \cite[p 287, Example 2]{Bo}). So, if $R$ is a Noetherian ring for which $\mathrm{Ass}(R)$ has embedded primes, then its total quotient ring $Q$ is a noetherian ring that is not Artinian. Hence, $Q$ does not have finite length as a $Q$-module. Hence, we have the following result.

\begin{thm}
If $R$ is Noetherian, there are finitely generated co-Hopfian $R$-modules that do not have finite length. \qed
\end{thm}

\begin{illus}
Let $R=k[x,y],$ $k$ a field. Let $I=(x^2,xy).$ The ideal $I$ has primary decomposition $(x) \cap (x,y)^2.$ Hence, the associated primes of $I$ are $P_1=(x)$ and $P_2=(x,y),$  and $P_1 \subset P_2.$ The ideal $I$ is not itself prime since $xy \in I$, but neither $x$ nor $y$ is. Let $Q$ be the total quotient ring of $R/I.$ Note that $Q$ has Krull dimension $1$. To see that $R/I$ does not have finite length, let $\bar{y}$ be the coset of $y$ in $R/I$ and consider the chain of ideals
$$(\bar{y}) \supsetneq (\bar{y}^2) \supsetneq (\bar{y}^3) \supsetneq \ldots .$$

This particular example can be generalized to any finite number of variables, yielding a $Q$ of any finite Krull dimension.

Let $R=k[x_1, \ldots,x_n].$ Let $I=(x_1^2,x_1x_2, \ldots,x_1x_n).$ Then $I$ has primary decomposition $(x_1)\cap(x_1, \ldots,x_n)^2.$ Hence, the associated primes of $I$ are $P_1=(x_1)$ and $P_2=(x_1,\ldots,x_n)$ with $P_1 \subset P_2.$ Again, $I$ is not prime: $x_1x_2 \in I$, but neither $x_1$ nor $x_2$ is. Let $Q$ be the total quotient ring of $R/I.$ Then $Q$ has Krull dimension $n-1$. To show $R/I$ does not have finite length, consider the chain of ideals generated by the integral powers of $\bar{x}_i^n$ for any $i \geq 2$.\qed
\end{illus}

\begin{remk}
Generally, co-Hopficity is independent of chain conditions: the abelian group $\mathbb{Z}(p^{\infty})$ has DCC but not ACC, while the $Q$-modules $Q$ of the preceding example have ACC but not DCC. Any finite abelian group has both chain conditions, while $\mathbb{Q},$ thought of as an abelian group, has neither.

For finitely generated modules, the first case listed above cannot occur. Let $R$ be commutative, $M$  a finitely generated $R$-module with $DCC$, and $I$ the annihilator of $M$. Then $R/I$ is an Artinian ring \cite[Lemma VIII.4.2]{Hung} and $M$ has $ACC$ as an $R/I$-module since $R/I$ is Noetherian. But then $M$ has $ACC$ as an $R$-module \cite[Prop 1.13]{S-V}. The last case cannot occur for a finitely generated module $M$ over a Noetherian $R$ since $M$ automatically has $ACC$.  \qed
\end{remk}

It is an entertaining exercise to characterize the finitely generated co-Hopfian $R$-modules which have finite length, $R$ Noetherian. The next theorem, which is elementary, lists several conditions equivalent to finite length in this context.

\begin{thm}
Let $R$ be commutative Noetherian, $M$ a finitely generated co-Hopfian $R$-module, and $I=\mathrm{ann}_R(M).$ The following are equivalent.

\begin{enumerate}
\item $M$ has finite length.
\item $R^{\prime}=R/I$ is Artinian.
\item $M$ has maximal orders.
\item $M$ is Artinian as an $R^{\prime}$-module.
\item $M$ is finitely cogenerated as an $R^{\prime}$-module.
\end{enumerate}
\end{thm}

\begin{proof}

$M$ is a faithful $R/I$-module. By Theorem 3.4, $R/I$ is semilocal and $\mathrm{Ass}_{R/I}(M)  \supseteq \Omega(R/I).$ For $M$ to have finite length, we must have $\mathrm{Spec}(R/I)=\Omega(R/I),$ i.e., $R/I$ must be Artinian. So, the equivalences $(1) \Leftrightarrow (2) \Leftrightarrow (4)$  are clear \cite[Cor. 2.17]{Eisenbud}.

$(1) \Rightarrow (3):$ Finite length implies DCC, implies maximal orders by Theorem 2.2.

$(3) \Rightarrow (5):$ By Theorem 2.3, $M$ is an essential extension of its socle and $M$ is a submodule of $\bigoplus _{\mathfrak{m} \in \Omega}(\oplus_{n_{\mathfrak{m}}}E(R/\mathfrak{m}))$, for cardinal numbers $n_{\mathfrak{m}}.$ But, the sum has only finitely many summands \cite[Thm. 4.7]{S-V}. So, the socle of $M$ is finitely generated and $M$ is finitely cogenerated.

$(5) \Rightarrow (1):$ $M$ is an essential extension of its socle, and the socle is finitely generated. So, $M$ has the form stated in Theorem 2.2.2. Hence, $M$ has DCC, and so finite length.
\end{proof}

\begin{remk}
The ring $Q$ from the example has annihilator $0$ as a $Q$-module. It does not have finite length. But $Q$ is not Artinian, so the example is consistent with the theorem.

Note that any one of the conditions of Theorem 4.9 implies that $M$ is co-Hopfian. The module $Q$ shows the converse is not true. \qed  \\
\end{remk}

A less obvious characterizaton of finitely generated co-Hopfian modules of finite length is that they are precisely those finitely generated $R$-modules $M$ which are {\em Cohen-Macaulay}. We were quite surprised by this connection, and stumbled on it somewhat by accident. There are a couple of things to keep in mind going forward. First, if $M$ is an $R$-module with annihilator $I$, then $M$ is co-Hopfian as an $R$-module if and only if it is co-Hopfian as an $R/I$-module (Theorem B). Secondly, $M$ has finite length as an $R$-module if and only if it has finite length as an $R/I$-module, since the submodule lattice is the same in either situation \cite[Prop 1.13]{S-V}. So, it suffices to examine what happens for {\em faithful} $R$-modules. 

It follows from Theorem 3.4 that if $M$ is finitely generated faithful and co-Hopfian, then the only $M$-regular elements of $R$ are units of $R$. So, no maximal ideal of {$R$ can contain an $M$-regular sequence \cite[p 123]{Mats}. Hence, the $\mathfrak{m}$-depth of $M$, $\mathrm{depth}(\mathfrak{m},M)$, is $0$ for each maximal ideal $\mathfrak{m}$ of $R$ \cite[Thm 16.7]{Mats}\cite[p.429]{Eisenbud}. The notion of depth is crucial in the definition of Cohen-Macaulay ($CM$) modules. 

If $R,\mathfrak{m}$ is local Noetherian, a finitely generated $R$-module $M$ is a $CM$-module if $\mathrm{depth}(\mathfrak{m},M)=\mathrm{dim}(R/I)$, $I=\mathrm{ann}(M).$ For general Noetherian $R$, the module $M$ is $CM$ if $M_{\mathfrak{m}}$ is $CM$ for each maximal ideal $\mathfrak{m} \in \mathrm{Supp}(M)$ \cite[Def 2.1.1]{B-H}. It turns out that finitely generated faithful co-Hopfian $CM$ modules are quite rare, as shown in the following pair of lemmas.

\begin{axe}
Let $R,\mathfrak{m}$ be local Noetherian and $M$ a finitely generated, faithful, $CM$ module. Then $M$ is co-Hopfian if and only $\mathrm{depth}(\mathfrak{m},M)=0$.
\end{axe}

\begin{proof}
$M$ is $CM$, so $\mathrm{depth}(\mathfrak{m},M)=\mathrm{dim}(R/\mathrm{ann}(M))=\mathrm{dim}(R)$, since $M$ is faithful (note $M$ is then a {\em maximal} $CM$ module (see \cite[Prop 21.9]{Eisenbud}).

$(\Rightarrow)$ If $M$ is co-Hopfian, then $\mathfrak{m}$ is an associated prime of $M$ (Theorem 3.4). So, $a \in R$ is $M$-regular if and only if $a$ is a unit. Hence, $\mathrm{depth}(\mathfrak{m},M)=0$ since $\mathfrak{m}$ cannot contain an $M$-regular sequence.

$(\Leftarrow)$ If $\mathrm{depth}(\mathfrak{m},M)=0$, then $R$ is Artinian and $M$ is co-Hopfian.
\end{proof} 

We have the expected generalization of Lemma 4.9.

\begin{axe}
Let $R$ be Noetherian and $M$ a finitely generated faithful $CM$ $R$-module. Then $M$ is co-Hopfian if and only if $\mathrm{depth}(\mathfrak{m}R_{\mathfrak{m}}, M_{\mathfrak{m}})=0$ for each maximal ideal $\mathfrak{m} \in \mathrm{Supp}(M)$.
\end{axe}

\begin{proof}
$(\Rightarrow)$ By Theorem 3.4, $R$ is semilocal and each of its maximal ideals is an associated prime of $M$, hence in $\mathrm{Supp}(M)$. Again, $a \in R$ is $M$-regular if and only if $a$ is a unit if and only if $a$ belongs to no maximal ideal of $R$. So $\mathrm{depth}(\mathfrak{m}, M)=0$ for each $\mathfrak{m} \in \mathrm{Supp}(M)$. But then $\mathrm{depth}(\mathfrak{m}R_{\mathfrak{m}}, M_{\mathfrak{m}})=0$ for each $\mathfrak{m}$  \cite[Lemma 18.1]{Eisenbud}. So, $M$ is $CM$.

$(\Leftarrow)$ Since $M$ is $CM$, $\mathrm{depth}(\mathfrak{m}R_{\mathfrak{m}}, M_{\mathfrak{m}})=0$ for each $\mathfrak{m}$ implies each $R_{\mathfrak{m}}$ is Artinian. By Lemma 4.9, the same facts imply each $M_{\mathfrak{m}}$ is co-Hopfian. Therefore, $M$ is co-Hopfian via local-global results \cite[Prop 3.9]{A-M}.
\end{proof}

We can now verify our main result.

\begin{thm}
Let $R$ be Noetherian and $M$ a finitely generated faithful co-Hopfian $R$-module. Then $M$ has finite length if and only if $M$ is $CM$.
\end{thm}

\begin{proof}
The hypotheses guarantee that $R$ is semilocal and that each of its maximal ideals is an associated prime of $M$ (Theorem 3.4). As before, $\mathrm{depth}(\mathfrak{m}, M)=0$ for each maximal ideal $\mathfrak{m}$.

$(\Rightarrow)$ If $M$  has finite length, then $R$ must be Artinian \cite[Lemma VIII.4.2]{Hung}. Hence, $\mathrm{dim}(R_{\mathfrak{m}})=0$ for each of the finitely many maximal ideals $\mathfrak{m}$ of $R$. But, $\mathrm{depth}(\mathfrak{m}R_{\mathfrak{m}}, M_{\mathfrak{m}})=0$ for each $\mathfrak{m}$. Hence, $M_{\mathfrak{m}}$ is Cohen-Macaulay for each $\mathfrak{m}$ whence $M$ is Cohen-Macaulay.

$(\Leftarrow)$ Since $M$ is $CM$, $M_{\mathfrak{m}}$ is $CM$ for each $\mathfrak{m}$. As seen in the proof of Lemma 4.10, each $M_{\mathfrak{m}}$ has $\mathrm{depth}(\mathfrak{m}R_{\mathfrak{m}}, M)= 0$ as an $R_{\mathfrak{m}}$-module. So, each $R_{\mathfrak{m}}$ is Artinian. But then $R$ itself is Artinian (any prime ideal properly contained in some $\mathfrak{m}$ yields a  prime properly contained in  $\mathfrak{m}R_{\mathfrak{m}}$). Hence, $M$ has finite length.
\end{proof}

 \section{Co-Hopfian submodules of $E(R/P)$}

Theorem 2.2 showed that if $P$ is a maximal ideal of $R$, then every submodule of $E(R/P)$ is co-Hopfian (for example, the abelian group $\mathbb{Z}(p^{\infty})$). What happens if $P$ is prime but not maximal? In this section, we use localization to identify the co-Hopfian  submodules of $E(R/P)$, actually of $\oplus_nE(R/P)$, for a given prime $P$. Later, we consider the more general case of an injective $E=\oplus_{i=1}^k(\oplus_{n_i}E(R/P_i))$ for a finite set of nonmaximal primes $P_1, \ldots, P_k$. With suitable restrictions on the $P_i$, we can determine the co-Hopfian submodules of such an $E$.

 In the proof of \cite[Thm 3.6]{Eben}, Matlis shows that if $P$ is a prime ideal of $R$, then $E(R/P)$ has the structure of an $R_P$-module. More generally, we have the following result \cite[Prop 5.6]{S-V}.

\begin{thm}
As an $R_P$-module, $E(R/P)$ is isomorphic to an $R_P$-injective envelope of $R_P/PR_P$, where $PR_P$ is the unique maximal ideal of $R_P$.  \qed
\end{thm}

\begin{remk}
If we consider$E(R/P)$ as a $R$-module, we will write $E_R(R/P)$ and if we consider it as an $R_P$-module write $E_{R_P}(R/P)$.
So, the isomorphism in the previous theorem would be written $E_{R_P}(R/P) \cong E_{R_P}(R_P/PR_P)$. The submodules of an injective envelope are to be considered as modules over the ring in the subscript. Likewise, all isomorphisms are module isomorphisms over the ring in the subscript. 
\end{remk} 

As mentioned above, all $(R_P$-)submodules of $E_{R_P}(R_P/PR_P)$ are co-Hopfian. We now have the answer to our question, and a bit more.

\begin{thm}
Let $P$ be a prime ideal of $R.$ If $n$ is a positive integer, then $M \subseteq \oplus_n E_R(R/P)$ is co-Hopfian if and only if $M$ is an $R_P$-module.
\end{thm}

\begin{proof}
Since the $0$ submodule of $\oplus_n E_R(R/P)$ is trivially co-Hopfian and an $R_P$-module, we assume $M\neq 0.$

$(\Rightarrow)$ Specialize Theorem 2.4 to the case in which $S=R\setminus P.$ 

$(\Leftarrow)$ If $M$ is an $R_P$-module, then $M$ is isomorphic to a submodule of the $R_P$-module $\oplus_n E_{R_P}(R_P/PR_P)$, which has $DCC$ (Theorem E). Hence $M$ is a co-Hopfian $R_P$-module, and so a co-Hopfian $R$-module (Theorem B).
\end{proof}

\begin{cor}
If $P$ is a nonmaximal prime of $R$, then the co-Hopfian submodules of $\oplus_n E_R(R/P)$ are the $R_P$-submodules of $\oplus_n E_{R_P}(R_P/PR_P)$ thought of as $R$-modules. \qed
\end{cor}

The module $E=E_R(R/P)$ has certain natural submodules. Since $E$ is $P$-primary, we may define, for each integer $n \geq 0$, an  $R$-submodule $E_n$ of $E$ by $E_n=\{x\in E:P^nx=0\}$. Note that $E_0=0$. It is clear that $E_n \subseteq E_{n+1}$ for all $n$. Moreover, $E=\cup E_n$ \cite[Thm 3.78]{Lam}. In addition, each $E_n$ is an $R_P$-module. We include the easy proof for completeness: $\mu_s$ is an automorphism of $E$ for each $s \in S=R\setminus P$; since $E_n$ is a fully invariant submodule of $E$, $\mu_s$ restricts to a monic endomorphism of $E_n$; if $y\in E_n$, then $y=\mu_s(x)$ for some $x \in E$; but $0=P^ny=P^n\mu_s(x)=\mu_s(P^nx)$. Since $\mu_s$ is monic, $P^nx=0$ and $x\in E_n$; thus $\mu_s$ is an automorphism of $E_n$ for each $s \in S$ and $E_n$ is an $R_P$-module. 

 The $E_n$ are analogues of the subgroups $\mathbb{Z}(p^n)$ of $\mathbb{Z}(p^{\infty})$ and they correspond to their counterparts $0:_{E_{R_P}(R_P/PR_P)}(PR_P)^n$ in $E_{R_P}(R_P/PR_P)$ \cite[p 102]{Lam}. So, they are co-Hopfian. By \cite[Cor 3.85]{Lam}, they are also Hopfian. We have the following.

\begin{cor}
For each nonnegative integer $n$, and for each nonzero prime ideal $P$, $E_n$ is a bi-Hopfian submodule of $E_R(R/P)$. \qed
\end{cor}

 The module $E_1$ plays an important role.
 
 \begin{thm}
 Let $P$ be a prime ideal of $R$. Then $E_1$ is the minimal co-Hopfian submodule of $E_R(R/P)$.
 \end{thm}
 
 \begin{proof}
 For a prime ideal $P$, $E_1=R_P/PR_P$, the fraction field of the integral domain $R/P$ \cite[Thm 3.4]{Eben}. Thus, if $P$ is maximal, then $E_1=R/P$, since $R_P \cong R$ in this situation. But then $E_1$ is a simple $R$-module. Hence, $E_1$ is the {\em unique} minimal nonzero module of $E$. To see this, let $M$ be a minimal nonzero module of $E$. Since $E$ is uniform, $M \cap E_1 \neq 0.$ But $M \cap E_1 \subseteq E_1$, which is simple. So, $M \cap E_1 = E_1.$ Hence, $M \supseteq E_1$. Since $M$ is minimal, $M=E_1.$ The result follows because $E_1$ is co-Hopfian.
 
 If $P$ is not maximal, then Theorem 4.2 implies that no proper nonzero $R$-submodule of $E_1$ is co-Hopfian.  By Theorem 5.3, $E_1$ must be the minimal nonzero co-Hopfian submodule of $E$.
 \end{proof}
 
 Note that if $P$ is not maximal, then $E_R(R/P)$ is guaranteed to have submodules that are not $R_P$-modules. In this case, $E_1=Q$, the fraction field of the domain $D=R/P$. But $Q \neq D$. The nonzero proper fractional ideals of $D$ are $D$-submodules of $Q$. But none are  $Q$-modules  since they are not co-Hopfian (Theorem 4.2 and introduction to Section 4). 
 
 The module $E_1$ is not the minimal Hopfian submodule of $E(R/P)$. Let $P$ be a nonmaximal ideal of $R$. Then $E_1$ is the fraction field of $R/P$, and properly contains it. But $R/P$ is a Hopfian $R/P$-module, so a Hofian $R$-module (Theorem B).    \\ 
 
 Let $\mathcal{P}$ be a set of primes of $R$. It would be nice to know when a submodule $M \subseteq E=\oplus_{P\in \mathcal{P}}(\oplus_{n_P}E(R/P)$ is co-Hopfian. We know $E$ is co-Hopfian. If we restrict the set $\mathcal{P}$ appropriately, we can determine the co-Hopfian submodules of $E$.  The following lemma suggests such a restriction. It is well-known and generalizes \cite[Prop 4.21]{S-V}.

\begin{axe}
Let $P,Q$ be prime ideals of $R$, neither of which contains the other. Let $M,N$ be, respectively, $P$-primary and $Q$-primary $R$-modules. Then $\mathrm{Hom}_R(M,N)=0$.
\end{axe}

\begin{proof}
Let $f:M \rightarrow N$ be a homomorphism and $x$ a nonzero element of $M$. Let $y=f(x)$ and suppose that $y \neq 0$. Then $P^nx=Q^my=0$ for some positive integers $n,m$. Hence, $P^ny=P^nf(x)=f(P^nx)=f(0)=0$. So, $P^n \subseteq \mathrm{ann}_R(y) \subseteq Q$. Hence, $P \subseteq Q$, a contradiction.
\end{proof}

Let $\mathcal{P}=\{P_i\}_{i \in I}$ be an {\em isolated} set of prime ideals of $R$ (as in \cite[\S 4.3]{S-V}): that is, no $P_i$ is contained in any $P_j$ for $i \neq j$. Let $M=  \oplus_{i \in I}M_i$ with each $M_i$ a $P_i$-primary $R$-module. Then any homomorphism $f:M \rightarrow M$ ``diagonalizes'' as $f=\oplus_{i \in I}f_{ii}$, with $f_{ii}  \in \mathrm{Hom}(M_i,M_i)$, since $\mathrm{Hom}(M_i,M_j)=0$ when $i \neq j$. We know from Theorem 2.2 that if $\mathcal{P}=\Omega(R)$, then every submodule of $E=\oplus_{P\in \mathcal{P}}(\oplus_{n_P}E(R/P)$ is co-Hopfian. What happens for other $\mathcal{P}$ can be determined via localization, similar to the treatment in Theorem 5.3.

Again let $E=\oplus_{P\in \mathcal{P}}(\oplus_{n_P}E(R/P)$ and $S=R\setminus (\cup_{P\in \mathcal{P}}P))$. Then $S$ is a multiplicative set \cite[Thm 2]{Kap}. Any $s\in S$ is nonzero and belongs to no $P \in \mathcal{P}$. Therefore, $\mu_s$ is an automorphism of each $E(R/P)$ \cite[Lemma 3.2]{Eben}\cite[Thm 18.4]{Mats}. Hence, each $E(R/P)$, and so each $\oplus_{n_P}E(R/P)$, is an $R_S$-module. By diagonalization, $\mu_s$ is an automorphism of $E$, and so $E$ is also an $R_S$-module. The ring $R_S$ has maximal ideals $PR_S$, $P \in \mathcal{P}$ \cite[Prop 3.11]{A-M}.

In a similar way, a co-Hopfian submodule $M$ of $E$ is an $R_S$-module. There is no harm in assuming that $E=E(M)$. In this case, $\mathcal{P}$ is the set of associated primes of $M$. Since any $s \in S$ belongs to no $P \in \mathcal{P}$, $\mu_s:M \rightarrow M$ is an injection. But $M$ co-Hopfian implies that $\mu_s$ is an automorphism of $M$.

Localization commutes with direct sums since tensor product does \cite[Thm 8.87]{Joe2}. So, $E_S=\oplus_{P\in \mathcal{P}}(\oplus_{n_P}[E_R(R/P)]_S)$. The module $E_R(R/P)$ is an $R_S$-module, and is injective as such \cite[Prop. 11.34]{Joe2}.  Now, $R/P \subseteq E_R(R/P)$, so we have  a chain of $R_S$-modules $(R/P)_S=R_S/PR_S \subseteq [E_R(R/P)]_S = E_{R_S}(R/P)$, the last equality because $E_R(R/P)$ is an $R_S$-module. So, we have the chain of $R$-modules $R/P \subseteq R_S/PR_S \subseteq E_R(R/P)$. Now $E_R(R/P)$ is essential over $R/P$, so it is essential over $R_S/PR_S$. Since it is an injective $R_S$-module, $E_{R_S}(R/P) \cong E_{R_S}(R_S/PR_S)$. Thus, $E_S=\oplus_{P\in \mathcal{P}}(\oplus_{n_P}E_{R_S}(R/P)))$.

If $M$ is an $R_S$-module, it can be considered as a submodule of $E_S$. Since the $\mathfrak{m}=PR_P$ are maximal ideals of $R_S$, and associated primes of $M_S=M$, $M$ has maximal orders and is the direct sum of its $\mathfrak{m}$-components.  But, each $\mathfrak{m}$-component of $M$ is co-Hopfian (Theorem 5.3). By diagonalization, $M$ is a co-Hopfian $R_S$-module, so a co-Hopfian $R$-module (Theorem B). 

\begin{thm}
Let $\mathcal{P}$ be an isolated set of prime ideals of $R$, not all of which are maximal. Let $S=R\setminus(\cup_{P\in \mathcal{P}}P)$. Let $E=\oplus_{P\in \mathcal{P}}(\oplus_{n_P}E(R/P))$, each $n_P\geq 0$ an integer. If $M$ is a submodule of $E$, then $M$ is co-Hopfian if and only if $M$ is an $R_S$-module. \qed \\
\end{thm}

There is a kind of ``primary decomposition'' underlying the theorem, and it is hinted at in the paragraph that precedes it. If $M$ is a submodule of the injective $E$ in the theorem, then for each $P \in \mathcal{P}$ we can define a subset $M(P)$ of $M$ as the set of $x \in M$ such that $P^nx=0$ for some nonnegative integer $n$. It is clear that $M(P)$ is a submodule of $E$. A few moments thought allows that $M(P)=M \cap(\oplus_{n_P}E_R(R/P))$, since $\oplus_{n_P}E_R(R/P)$ is precisely the set of $x \in E$ which are annihilated by a power of $P$. Consequently, $M=\sum_{P\in \mathcal{P}}M(P)$. Of course, we want the sum to be direct. To this end, suppose there is a nonzero $x$ in the intersection $M(P)\cap (\sum_{Q\neq P} M(Q))$. Since $M(P)$ is $P$-primary, so is the cyclic submodule $Rx$ generated by $x$. Now $P^n \subseteq \mathrm{ann}_R(x)=\mathrm{ann}_R(Rx)$, so $P \subseteq \mathfrak{r}(\mathrm{ann}_R(x))\subseteq P$, since $P$ is the only prime containing $\mathrm{ann}_R(x)$, $\mathfrak{r}$ denoting radical \cite[Thm 3.1]{Eisenbud}. Thus, $\mathrm{ann}_R(x)$ is a $P$-primary ideal.

But $x$ also belongs to $\sum_{Q\neq P}M(Q)$. So, $x=\sum_{i=1}^k x_i$ with $x_i \in Q_i$. Now, $\mathrm{ann}_R(x)=\cap_{i=1}^k\mathrm{ann}_R(x_i)$. Taking radical of both sides yields $P=\cap_{i=1}^kQ_i$. Therefore, $P\supseteq Q_i$ for some $i$ \cite[Prop 1.11]{A-M}, contradicting the fact that $\mathcal{P}$ is an isolated set of primes. Since $P$ was chosen arbitrarily, $M=\oplus_{P\in \mathcal{P}}M(P)$. Thus, we have a result analogous to Theorem F(3).

\begin{thm}
Let $\mathcal{P}$ be an isolated set of prime ideals of $R$, and $E$ the co-Hopfian injective $R$-module of Theorem 5.8. If $M\subseteq E$, then $M=\oplus_{P\in \mathcal{P}}M(P)$, with $M(P)= M \cap(\oplus_{n_P}E_R(R/P))$. \qed
\end{thm}

The result in Theorem 5.8 is similar to the ``cutting property'' of quasi-injective modules pointed out by Lam \cite[Prop 6.79]{Lam}.\\

 Let $E$, $M$, and $M(P)$ be as in Theorem 5.9. Let $E(P)=\oplus_{n_P}E_R(R/P)$. If $M$ is an $R_S$-module, then $M \cong M_S = \oplus_{P\in \mathcal{P}}M(P)_S$. Since localization is exact \cite[Prop 3.3]{A-M}, $M(P)_S \subseteq E(P)_S$ and $M(P)_P \subseteq E(P)_P$. For any prime $P$, $S \subseteq R \setminus P$, so $R_P$ is a localization of $R_S$ \cite[Thm 4.3 Cor 3]{Mats}. But then $M(P)_P$ is a localization of $M(P)_S$ (stacks.math.columbia.edu/tag/00CM). By Theorem B, $M(P)_P$ and $M(P)_S$ have the same endomorphism ring. But $M(P)_P$ is an $R_P$-submodule of $E(P)_P$, so is co-Hopfian by Theorem 5.3. Hence, $M(P)_S$ is co-Hopfian (Theorem B, again). But then $M_S \cong M$ is co-Hopfian by Theorem A and Lemma 5.7. So, Theorem 5.9 provides an alternate approach to proving the  $M=M_S$ implies $M$ co-Hopfian portion of Theorem 5.8.

\section{Hopfian submodules of $E(R/P)$}.

The results in the Hopfian case differ from those in the co-Hopfian case, even for $E(R/P)$. If $R$ is Artinian, then the indecomposable injective $E(R/\mathfrak{m})$ is Hopfian for each maximal ideal $\mathfrak{m}$  of $R$ \cite[Thm 3.5]{Me3}.  More generally, we showed in Section 2 of that article that if $R$ is a Noetherian domain, then $E(R/P)$ is Hopfian if and only if $P$ is the zero ideal, and, more generally, if $R$ is Noetherian then $E(R/P)$ is Hopfian if $P$ is a minimal prime of $R$. All minimal primes of a Noetherian $R$ are associated primes of $R$. In fact, the associated primes of $R$ provide the answer in general. 

\begin{thm}
Let $R$ be Noetherian and $n$ a positive integer. Then $\oplus_nE(R/P)$ is Hopfian if and only if $P$ is an associated prime of $R$. 
\end{thm}

\begin{proof}
The proof requires some care. Ultimately, we will reduce to the following specific case. Let $R, \mathfrak{m}$ be a local, Noetherian, complete quoring, and let $E=E(R/\mathfrak{m})$. The module $E$ is the minimal injective cogenerator of the category of $R$-modules. Hence, the contravariant functor $D(-)=\mathrm{Hom}_R(-,E)$ is {\em faithfully exact}, that is, a given sequence of $R$-modules $0 \rightarrow A \rightarrow B \rightarrow C \rightarrow 0$ is exact if and only if the sequence $0 \rightarrow D(C) \rightarrow D(B) \rightarrow D(A) \rightarrow 0$ is exact. $D$  is the functor that assigns to any $R$-module $M$ its {\em Matlis dual}. 

If $f:\oplus_nE \rightarrow \oplus_nE$ is a surjection, then $\oplus_nE \stackrel{f}{\rightarrow} \oplus_nE \rightarrow 0$ is exact. Conveniently, $D(E)=\hat{R}$, the $\mathfrak{m}$-adic completion of $R$ \cite[Thm 5.14 and Cor 1,2]{S-V}. Since $R$ is complete, $D(E)=R$. So, $D(\oplus_nE)=R^n$. Hence, $0 \rightarrow R^n \stackrel{f*}{\rightarrow} R^n$ is exact. Since $R$ is a quoring, $R^n$ is co-Hopfian \cite[Thm 2.5]{Me}. So, the sequence $0 \rightarrow R^n \stackrel{f*}{\rightarrow} R^n \rightarrow 0$ is exact. Hence, by the faithful exactness of $D$, the map $f:\oplus_nE \rightarrow \oplus_nE$ is an isomorphism and $\oplus_nE$ is Hopfian.\\

 First, we observe that $\oplus_nE(R/P)$ is Hopfian only if $P$ is an associated prime of $R$.

\begin{axe}
Let $R$ be Noetherian. If $P \not\in \mathrm{Ass}(R)$, then $\oplus_nE(R/P)$ is not Hopfian.
\end{axe}

\begin{proof}
The zero divisors of $R$ are the elements in the union of the associated primes of $R$. So, if $P$ is not an associated prime, prime avoidance yields an element $r \in P$ that is a regular element of $R$. Since $\oplus_nE(R/P)$ is divisible (it's injective), $r(\oplus_nE(R/P))=\oplus_nE(R/P)$. But $\mu_r$ is not an injection, since $r \in P$ and $P= \mathrm{ann}(R/P)$. 
\end{proof}

The next lemma tells us when a local ring is a quoring.

\begin{axe}
Let $R, \mathfrak{m}$ be local Noetherian. The conditions (a) $R$ is a quoring; (b) $R$ is co-Hopfian; and (c) $\mathfrak{m} \in \mathrm{Ass}(R)$ are equivalent.
\end{axe}

\begin{proof}
The equivalence of $(a)$ and $(b)$ was noted in Section 1. The equivalence of $(a)$ and $(c)$ was proved in Theorem 3.6.
\end{proof}

Now, reduce to the local case. The $M=R$ case of \cite[Thm 6.2 Cor]{Mats} tells us that  if $R$ is Noetherian and $P \in \mathrm{Spec}(R)$, then $P \in \mathrm{Ass}_R(R) \Leftrightarrow PR_P \in \mathrm{Ass}_{R_P}(R_P).$ In particular, if $P$ is an associated prime of $R$, then $R_P$ is a local quoring (Lemma 6.3). We know $E_R(R/P)$ is an $R_P$-module isomorphic to $E_{R_P}(R_P/PR_P)$ (Theorem 5.1). The canonical map $R \rightarrow R_P$ is a ring epimorphism, so $E_R(R/P)$ is Hopfian as an $R$-module if and only if it is Hopfian as an $R_P$-module (Theorem B). So, we may assume that $R, \mathfrak{m}$ is a  local quoring and deal with $E=E_R(R/\mathfrak{m})$. 

In passing to the completion $\hat{R},\hat{\mathfrak{m}}$, we need to know two things. First, is $\hat{\mathfrak{m}}$ an associated prime of $\hat{R}$? Since $\mathfrak{m}$ is an associated prime of $R$, we have an exact sequence $0 \rightarrow R/\mathfrak{m} \rightarrow R$. Since $R$ is local Noetherian, $\hat{R}$ is a faithfully flat $R$-module \cite[Thm 8.14]{Mats}. So, tensoring the previous sequence over $R$ with $\hat{R}$ yields the exact sequence $0 \rightarrow \hat{R}/\hat{\mathfrak{m}} \rightarrow \hat{R}$, since $\widehat{R/\mathfrak{m}}=\hat{R}/\hat{\mathfrak{m}}$ \cite[Thm 8.1]{Mats}. So, $\hat{\mathfrak{m}}$ is an associated prime of $\hat{R}$. Hence, $\hat{R}, \hat{\mathfrak{m}}$ is a complete, local, Noetherian quoring.

Secondly, we need to know how the Hopfian property behaves under the change of ring from $R$ to $\hat{R}$. The $R$-module $E$ is an $\hat{R}$-module, and  $\mathrm{Hom}_R(E,E)=\mathrm{Hom}_{\hat{R}}(E,E)=\hat{R}$ \cite[Thm 18.6]{Mats}. So, $E$ is Hopfian as an $R$-module if and only if it is Hopfian as an $\hat{R}$-module. But $\oplus_nE$ is also an $\hat{R}$-module.  Its endomorphism ring is the set of $n \times n$ matrices with entries that are maps $E_R(R/\mathfrak{m})\rightarrow E_R(R/\mathfrak{m})$ \cite[Prop 1.24]{S-V}, but, as mentioned above, these maps are scalar multiplications by elements of $\hat{R}$, and so are both $R$-linear and $\hat{R}$-linear. In other words, $\mathrm{Hom}_R(\oplus_n E_R(R/\mathfrak{m}),\oplus_n E_R(R/\mathfrak{m}))=\mathrm{Hom}_{\hat{R}}(\oplus_n E_R(R/\mathfrak{m}),\oplus_n E_R(R/\mathfrak{m}))$. Hence, the sum $\oplus_nE$ is Hopfian as an $R$-module if and only if it is Hopfian as an $\hat{R}$-module. The proof of our initial case now applies to show $\oplus_nE$ is Hopfian as $\hat{R}$-module. Lemma 5.8 established the converse, so Theorem 6.1 is proved.
\end{proof}

Consequently, Hopfian injective modules of the form $\oplus_nE(R/P)$ are rare if $R$ is a Noetherian ring with infinitely many prime ideals.  Since $R$ has only finitely many associated primes and $n$ is a positive integer, there are only countably many Hopfian injectives $\oplus_nE(R/P)$ (and these are the only bi-Hopfian injective $R$-modules). But, there are at least an uncountable number of injectives of that form that are not  Hopfian. \\

Suppose we have a set of prime ideals $\mathcal{P}$ in an arbitrary commutative Noetherian ring $R$. Let $E=\oplus_{\mathcal{P}}(\oplus_{n_P}E_R(R/P))$, with each $n_P$ a nonegative integer. We know that $E$ is co-Hopfian (Theorem C). If we want $E$ to be Hopfian, we will need at least to restrict  our primes to associated primes of $R$. If we restrict to an isolated set of associated primes of $R$, we can say something meaningful.

There are commutative Noetherian rings $R$ that have an isolated set of associated primes, for instance Cohen-Macaulay rings and reduced rings (e.g., an integral domain). For other $R$ we will simply choose some isolated subset of the associated primes. Since $R$ is Noetherian, all such subsets are necessarily finite. 

\begin{thm}
Let $\mathcal{P}=\{P_1, \ldots, P_k\}$ be an isolated set of associated primes of $R$. Let $E=\oplus_{P\in \mathcal{P}}(\oplus_{n_i}E_R(R/P_i) )$ with each $n_i$ a nonnegative  integer. Then $E$ is Hopfian.
\end{thm}

\begin{proof}
For convenience, let $M_i=\oplus_{n_i}E_R(R/P_i)$. By Theorem 6.1, each $M_i$ is Hopfian. Let $f:E \rightarrow E$ be a surjection. By Lemma 5.7, $f$ diagonalizes as $f= \oplus_{i=1}^kf_{ii}$ with $f_{ii}:M_i\rightarrow M_i$ (also a surjection). Since $M_i$ is Hopfian, $f_{ii}$ is an isomorphism. So, $f$ is, too, and $E$ is Hopfian.
\end{proof}

Let $E$ be the injective module in Theorem 6.4. What nonzero submodules $M$ of $E$ are Hopfian? Since the set $\mathcal{P}$ of primes is isolated, it suffices to answer the question for $E= \oplus_{n}E_R(R/P)$, $P$ a specific prime in $\mathcal{P}$ (Lemma 5.7). Should $R$ be hereditary, the answer is simple: every submodule of $E$ is Hopfian (0 and $E$ are automatic).

\begin{thm}
Let $R$ be commutative, Noetherian, and hereditary, $P$ an associated prime of $R$. Every submodule $M \subseteq \oplus_nE(R/P)$ is Hopfian.
\end{thm}

\begin{proof}
There is no harm in assuming that $E=E(M)$. Let $\varphi:M \rightarrow M$ be a surjection. Then $\varphi$ lifts to $\bar{\varphi}:E \rightarrow E$. We have $M=\bar{\varphi}(M)\subseteq \bar{\varphi}(E) \subseteq E$. Since $R$ is hereditary, $\bar{\varphi}(E)$ is an injective module. Hence, $\bar{\varphi}(E)=E$ since $E=E(M)$. So, as pointd out in the proof of Theorem 6.1, $\bar{\varphi} \in \mathrm{Mat}(n,\hat{R})$, $\hat{R}$ the $P$-adic completion of $R$. But $\mathrm{Mat}(n,\hat{R})$ is also the endomorphism ring of the free $\hat{R}$-module $\oplus_n\hat{R}$. Hence, $\bar{\varphi}$ is an isomorphism \cite[Thm III.8.2.1]{Bourbaki}. Thus, $\varphi:M \rightarrow M$ is an injection, and $M$ is Hopfian.
\end{proof}

The proof is nice, but it obscures the fact that the module $E$ in the theorem is simply a finite direct sum of fields. The ring $R$ of the theorem is a finite direct sum of rings that are Dedekind domains, some of which may be fields (McConnell and Robson \cite[\S 5.4.1]{M-R} mention it, and refer to \cite[Thm 168]{Kap}; R. van Dobben de Bruyn gives a nice proof at https://math.stackexchange.com.net/q/311603 that connects the two). 

Write $R=\oplus_{j=1}^nR_j$. It is well-known that an ideal $I$ of $R$ has the form $I=\oplus_{j=1}^nI_j$ with $I_j$ an ideal of $R_j$.  A prime ideal $P$ of $R$ has the same form with $I_k=P_k$ a prime ideal of $R_k$, and $I_j=R_j$ for $j\neq k$. So the quotient $R/P$ is nothing but the quotient $R_k/P_k$. If $P_k$ is the 0-ideal of $R_k$, then the prime ideal $P$ just defined will be an associated prime of $R$, annihilating any element $r=(r_1, \ldots ,r_n)$ of for which $r_k$ is a nonzero element of $R_k$, and $r_j=0$ for $j\neq k$. Every associated prime of $R$ arises in this fashion; there are exactly $n$ of them. Clearly, this set of associated primes is isolated. Note that even if $R_j=R_k$ for some $j,k$ and the primes $P_j=P_k=0$, the corresponding primes $P$ and $P^{\prime}$ of $R$ are nonetheless isolated.

So, $E_R(R/P) = E_{R_k}(R_k/0)$ is the quotient field $Q_k$ of the Dedekind domain $R_k$, or is the field $Q_k=R_k$ if $R_k$ is a field. Therefore, $E$ is the $n$-dimensional $Q_k$-vector space $Q_k^n$. Theorem 4.4 finishes the argument for $\oplus_nE(R/P)$. Lemma 5.7 allows us to say that the general $E=\oplus_{P\in \mathcal{P}}(\oplus_{n_k}E_R(R/P) )$, with each $n_k$ a nonnegative  integer, is the direct sum $\oplus_{P\in \mathcal{P}}Q_k^{n_k}$, where $\mathcal{P}=\{P_1, \ldots , P_n\}$ is the set of associated primes of $R$. By Theorem 5.9, the submodules of $E$ are the $\oplus_{P\in \mathcal{P}}Q_k^{m_k}$ with $m_k \leq n_k$ for all $k$.

Theorem 6.5 applies to the specific case in which $R$ is a field, or a Dedekind domain that is not a field.\\

Things get more interesting if we allow our set $\mathcal{P}$ to contain primes that are not associated primes of $R$. Virtually any set $\mathcal{P}$ of primes of $R$ is isolated. The only restriction necessary is we cannot have two prime ideals one of which uses the 0-ideal and the other a maximal ideal of the same ring $R_j$ if $R_j$ is a Dedekind domain.

Given our direct sum decomposition of $R$, a given isolated set $\mathcal{P}$ of primes of $R$ can be partitioned into sets $\mathcal{P}_j$ of those primes $P$ of $R$ that are built using an isolated set of primes of the summand  $R_j$ of $R$. Therefore, the general co-Hopfian injective $R$-module $E$ can be expressed as $E=\oplus_{j=1}^n(\oplus_{P\in \mathcal{P}_j}(\oplus_{n_P}E(R_j/P)))$, each $n_P$ a nonnegative integer. Note that $R_j/P$ is really $R_j/P_j$ for an appropriate prime of $R_j$. By the decomposition result Theorem 5.9, it suffices to determine which submodules of $E=\oplus_{n_P}E(R_j/P)$ are Hopfian for the various $P \in \mathcal{P}$. If $P=0$, we already know the answer from Theorem 4.4. It remains to consider the case in which $R_j$ is a Dedekind domain and $P$ a maximal ideal of $R_j$. Then $E$ is not Hopfian since $E(R_j/P)$ is not. If any submodule $M$ of $E$ contains an injective submodule, then $M$ is not Hopfian (the injective submodule is a non-Hopfian summand (Theorem A)). Any $M$ that does not have have an injective submodule is Hopfian. The proof of Theorem 6.5 works. We can sum things up as follows.

\begin{thm}
Let $\mathcal{P}$ be an isolated set of primes of the commutative, Noetherian, hereditary ring $R$. Let $E=\oplus_{j=1}^n(\oplus_{P\in \mathcal{P}_j}(\oplus_{n_P}E(R_j/P)))$, each $n_P$ a nonegative integer. A submodule $M\subseteq E$ is Hopfian if and only if its only indecomposable injective submodules are of the form $E(R/P)$ for $P$ an associated prime of $R$. These $M$ are, in fact, bi-Hopfian. \qed
\end{thm}

\section{$E(R/\mathfrak{m})$ with few co-Hopfian modules}

Any indecomposable injective $E(R/P)$ for $P$ a prime ideal of $R$ always has the $E_n$ as submodules. We can inquire as to when these are the only submodule of $E(R/P)$. It is a straightforward exercise in commutative algebra to determine for which rings $R$, always commutative Noetherian, this is true. It is easy to extend the analysis to decide for which commutative Noetherian $R$ the $E_n$ are the only submodules of $E(R/P)$ for each prime ideal $P$ of $R$.

 By Corollary 5.4, we can reduce the problem of determing the co-Hopfian submodules of $E_R(R/P)$  to determining the submodules of $E_{R_P}(R_P/PR_P) $. But, $R_P,PR_P$ is a local ring, so we may as well work with $E=E_R(R/\mathfrak{m})$ for $R, \mathfrak{m}$ local. We always have the filtration
\begin{eqnarray}
0=E_0 \subseteq E_1 \subseteq E_2 \subseteq \ldots \subseteq E_n \subseteq \ldots \subseteq E= \cup E_n
\end{eqnarray}
of $E$ by the submodules $E_n.$  From the proof of Theorem 5.6, we know that $E_1(=R/\mathfrak{m})$ is a simple $R$-module and the minimal submodule of $E$. 

There is a companion filtration of $R$ given by the annihilators of the $E_n$
\begin{eqnarray}
0=\cap \mathfrak{m}^n \subseteq \ldots \subseteq \mathfrak{m}^n \subseteq \ldots \mathfrak{m}^2 \subseteq \mathfrak{m} \subseteq R.
\end{eqnarray}

Since $R$ is local Noetherian, there are two possibilities for the filtration $(2)$: the $\mathfrak{m}^n$ are all distinct; or, $\mathfrak{m}^n = 0$ for some $n$, in which case $R$ is local Artinian \cite[Prop 8.6]{A-M}. In the latter case, there are clearly only finitely many $E_n$. In the former, the $E_n$ are all distinct. These observations are consequences of a close relationship between the two filtrations.

Let $k=R/\mathfrak{m}$. It is well-known that the filtration factors $E_n/E_{n-1}$ of $(1)$ are finite dimensional $k$-vector spaces \cite[Thm 3.4, 3.9]{Eben} (or \cite[p 102]{Lam}). In fact, for each positive integer $n$, $E_n/E_{n-1} \cong \mathfrak{m}^{n-1}/\mathfrak{m}^n$ as a $k$-vector space \cite[Thm 3.82]{Lam}. This observation justifies the following lemma.

\begin{axe}
Let $R,\mathfrak{m}$ be local Noetherian, $E=E(R/\mathfrak{m})$, and the $E_n$ as already defined. There is an $R$-module lying properly between $E_{n-1}$ and $E_n$ if and only if there is an ideal of $R$ lying properly between $\mathfrak{m}^n$ and $\mathfrak{m}^{n-1}$.
\end{axe}

Theorem 5.6 is the $n=1$ case of this lemma. Moreover, the lemma shows we can refine (in the usual sense) one of the filtrations (1), (2) if and only if we can refine the other. So, the $E_n$, $n \geq 1$, are the only proper nonzero submodules of $E_R(R/\mathfrak{m})$ if and only if the $\mathfrak{m}^n$ are the only proper nonzero ideals of $R$ (with appropriate interpretation if $R$ is Artinian). For what local Noetherian rings $R$ does this happen? 

The filtration $(2)$ cannot be refined if and only if there is no ideal of $R$ lying properly between $\mathfrak{m}^2$ and $\mathfrak{m}$. Necessity is clear. Sufficiency is an exercise in basic commutative algebra. If there is no ideal properly between $\mathfrak{m}^2$ and $\mathfrak{m}$, then $\mathrm{dim}_k\, \mathfrak{m}/\mathfrak{m}^2$, the {\em embedding dimension} of $R$, is either $0$ or $1$. If  $\mathfrak{m}/\mathfrak{m}^2=0$, then $\mathfrak{m}=0$ (Nakayama) and $R$ is a field. In this case, $E_0=0$ and $E_1=R$ are the only $R$ submodules of $E$. If $R$ is not a field, then $\mathfrak{m} \neq \mathfrak{m}^2$ (also Nakayama), and the embedding dimension of $R$ is $1$. Hence, $\mathfrak{m}$  is a principal ideal \cite[Thm 2.3]{Mats}. But then every proper nonzero ideal of $R$ is principal (adapt the proof of \cite[Prop 8.8]{A-M}, replacing $\mathfrak{m}$ is nilpotent by $\cap \mathfrak{m}^n=0$, the latter follows by \cite[Cor 10.19]{A-M}; or see \cite[Prop 8.3]{Reid} and the remark following it).

There are now two possibilities to consider since, in general, the Krull dimension of $R$ is no larger than the embedding dimension \cite[p 104]{Mats}. If $0$ is a prime ideal, then $R$ is $1$-dimensional and so is a $DVR$ \cite[Prop 9.2]{A-M}. If $0$ is not prime, $R$ is Artinian, $\mathfrak{m}$ is nilpotent, and $R$ is what Zariski and Samuel call a {\em special} $PIR$ \cite[p 245]{Z-S}. In the first case case, all containments in (1) are strict and the filtration has countable length. In the second, as alluded to earlier, all containments are strict, but the filtration is finite, with $E_n=E$, $n$ the {\em index of nilpotency} of $\mathfrak{m}$ (the least positive integer $n$ for which $\mathfrak{m}^n=0$). So, we have the following result.

\begin{thm}
Let $(R,\mathfrak{m})$ be local. The modules $E_n$ are the only submodules of $E=E(R/\mathfrak{m})$ if and only if $R$ is a $PIR$ of Krull dimension at most $1$. \qed
\end{thm}

Given a Noetherian ring $R$, we can ask when, for each prime $P$ of $R$, the $E_n$ are the only submodules of $E=E_{R_P}(R_P/PR_P)$. Theorem 7.2 implies that each $R_P$ must be either a $DVR$ or a special $PIR$. What rings $R$ have this property? We argue, in part, as suggested in the answer that  user math54321 provided in \cite{Me5}. 

Any Noetherian ring $R$ decomposes as a finite direct sum  $R=\oplus_{j=1}^n R_j$ of irreducible rings \cite[Prop 22.2]{Lam3}, \cite[Thm 7.9]{A-F} (finite direct sum and finite direct product are the same).

The rings singled out in Theorem 7.2 have the property that every nonzero prime is maximal. As we noted in Section 6, if $R=\oplus_{j=1}^nR_j$, then a prime ideal of $R$ has the form $P=\oplus_{j=1}^nP_j$ with $P_k$ a prime ideal in $R_k$ for some $k$, $1 \leq k \leq n$, and $P_j=R_j$ for $j \neq k$.  As usual, $R\setminus P$ is a multiplicative set in $R$ and is equal to $\oplus_{j=1}^nS_j$ where $S_k=R_k\setminus P_k$ and $S_j=R_j$ for $j \neq k$. A straightforward argument shows that $R_P=\oplus_{i=1}^n(R_j)_{S_j}=(R_k)_{P_k}$, since $0 \in S_j$ when $j \neq k$. So, localizations of $R$ at a prime corresponds to  the localization of some $R_j$ at one of its primes. Since each $R_P$ has dimension at most $1$, each $(R_j)_{P_j}$ has dimension at most $1$. In particular, each maximal ideal of each $R_j$ has height at most $1$. Thus, $R$ itself has dimension no greater than $1$.

If $R_j$ has dimension $0$, then it is Artinian, so a finite direct product of local Artinian rings. But $R_j$ is irreducible, so $R_j$ is local Artinian. Let $\mathfrak{m}$ be its maximal ideal. Then $(R_j)_{\mathfrak{m}}=R_j$ since the elements of $R_j$ not in $\mathfrak{m}$ are already units. Hence, $R_j$ is a $PIR$, since $(R_j)_{\mathfrak{m}}$ is. In fact, $R_j$ is a special $PIR$.

For notational simplicity, let $R=R_j$. If $R$ has dimension 1, then $R$ has a maximal ideal $\mathfrak{m}$ of height 1. Then the local ring $R_{\mathfrak{m}}$ also has dimension 1. By assumption, the maximal ideal $\mathfrak{m}R_{\mathfrak{m}}=\mathfrak{n}$ of $R_{\mathfrak{m}}$ is principal. Say $\mathfrak{n}=(t)$. But then every nonzero proper ideal $I$ of $R_{\mathfrak{m}}$ is a power of $\mathfrak{n}$. We adapt the proof of  \cite[Prop 8.8]{A-M} . It is clear that $I \subseteq \mathfrak{n}$. By Krull's Intersection Theorem \cite[Cor 5.4]{Eisenbud} \cite[Cor 10.19]{A-M}, $\bigcap_{j=1}^{\infty} \mathfrak{n}^j=0$. So, there is a $j$ such that $I \subseteq \mathfrak{n}^j$, but $I \not\subseteq \mathfrak{n}^{j+1}$. Choose $x \in I \setminus \mathfrak{n}^{j+1}$. Then $x=ut^j$ for some unit $u$ of $R_{\mathfrak{m}}$ (or $x \in \mathfrak{n}^{j+1}$). But then $t^n \in I$, so that $(t^j)=\mathfrak{n}^j \subseteq I$, and $I = \mathfrak{n}^j$.

 If $P$ is a nonzero prime ideal of $R_{\mathfrak{m}}$, then $P=\mathfrak{n}^k$ for some $k$. Taking the radical of both sides yields $P=\mathfrak{n}$, so that $R_{\mathfrak{m}}$ has $\mathfrak{n}$ as its unique nonzero prime. Since $R_{\mathfrak{m}}$ has dimension 1, 0 must be a prime ideal of $R_{\mathfrak{m}}$. So, $R_{\mathfrak{m}}$ is a domain, in fact, a $DVR$. 
 
 If all maximal ideals of $R$ have height 1, we are done \cite[Thm 9.3]{A-M}. Generally, there is no guarantee that all maximal ideals in a ring have the same height. For example, if $A=k[x,y,z]$, $k$ a field, let $P=(x)$ and $Q=(y,z)$. Let $S=A\setminus (P \cup Q)$. Then $A_S$ has maximal ideals $P_S$ and $Q_S$, the first having height 1, the second, height 2 \cite[Prop II.3.5.17]{Bourbaki}. The ring $A_S$ is a domain, so it is irreducible since it has no nontrivial idempotents \cite[Cor 7.7]{A-F}(who use the term {\em indecomposable}).
 
 In our situation, the rings $R$ are irreducible of dimension 1. In this case, all maximal ideals will have height 1. In any event, the only possible heights for a prime ideal of $R$ are 1 or 0.  We need only show the second possibility does not occur for a maximal ideal. To do so, we exploit some features of the prime spectrum $\mathrm{Spec}(R)$ of the ring $R$ (the basic facts are set forth in Exercises 15-26 on pages 12-15 of \cite{A-M}). Recall that $R$ represents one of the rings $R_j$.
 
 We first observe that since $R$ is irreducible, its spectrum is connected (Exercise 22). So, $\mathrm{Spec}(R)$ is not the union of two nonempty open disjoint sets. We will show that a maximal ideal of height 0 leads to a contradiction of this fact.
 
 Since $R$ is Noetherian, it has only finitely many minimal primes $P_1, \ldots , P_n$. Note that there are no containments between the various $P_j$, and that all these primes have height 0. Let $P_1, \ldots,P_k$ be those that are not maximal, and $P_{k+1}, \ldots,P_n$ those that are maximal. Let $I=\bigcap_{j=1}^k P_j$ and $J=\bigcap_{j=k+1}^nP_j$. The set $V(I)$ consists of all prime ideals that contain $I$. So, $V(I)$ contains $P_1, \ldots,P_k$ as well as all height 1 primes. However, $V(J)$ contains only the height 0 primes of $R$. Both $V(I)$ and $V(J)$ are closed sets in the Zariski topology of  $\mathrm{Spec}(R)$ (Exercise 15). Moreover, $V(I)$ and $V(J)$ are disjoint. In addition, they are complements of each other in  $\mathrm{Spec}(R)$. Hence, both are also open sets in  $\mathrm{Spec}(R)$. The set $V(I)$ is always nonempty. So, if there is a height 0 maximal ideal in $R$ we contradict  $\mathrm{Spec}(R)$ connected. Thus, all maximal ideals of $R$ have height 1. 
 
 All minimal ideals have height 0. Let $P$ be a minimal prime. Then $P$ is properly contained in some maximal ideal $\mathfrak{m}$ of $R$. So, in $R_{\mathfrak{m}}$, $P_{\mathfrak{m}}$ is distinct from $\mathfrak{m}_{\mathfrak{m}}$. Since $R_{\mathfrak{m}}$ is a $DVR$, $P_{\mathfrak{m}}$ must be the ideal 0. If $P$ is not contained in the maximal ideal $\mathfrak{n}$, then $P_{\mathfrak{n}}=0$ since $P$ has nontrivial intersection with $S=R \setminus \mathfrak{n}$. Thus, $P$ localizes to 0 for every maximal ideal in $\Omega$, the maximal spectrum of $R$. Hence, $P=0$ \cite[Prop 3.8]{A-M}, and $R$ is a domain. Since the dimension of $R$ is 1, $R$ is a Dedekind domain \cite[Thm 9.3]{A-M}. Thus, we have the following result.

\begin{thm}
Let $R$ be commutative Noetherian. Then for each prime ideal $P$ of $R$, the $E_n$ are the only submodules of $E=E_R(R/P)$ if and only if $R$ is a finite direct sum of Dedekind domains and special $PIR$s. \qed
\end{thm}

\begin{remk}
If the Dedekind domains in the decomposition of $R$ are local, they are $DVRs$, in which case $R$ itself is a $PIR$ (\cite[Thm IV.13.33]{Z-S}, which offers a decomposition characterizing $PIRs$ (recall that finite direct sums are isomorphic to finite direct products). \qed
\end{remk}

We have already pointed out that the the cyclic groups $\mathbb{Z}(p^n)$ are the modules $E_n$ for the indecomposable injective $\mathbb{Z}(p^{\infty})=E(\mathbb{Z}/p\mathbb{Z})$. This feature carries over to the $E_n$ as specified in Theorem 7.3. So, let $R$ be a finite direct sum of Dedekind domains and special $PIR$s, and let $P$ be a prime ideal of $R$. As usual, the localization of $R$ at $P$ is either a $DVR$ or a special $PIR$. We begin with the case of a local ring $R,\mathfrak{m}$ that is a $DVR$. 

We know that the submodules $E_n$, $n \geq 1$ of $E=E_R(R/\mathfrak{m})$ are distinct and nonzero. It is clear that $E_1$ is cyclic, since $E_1 \cong R/\mathfrak{m}$. For $n \geq 2$, the result will follow from the isomorphism $E_n/E_{n-1} \cong \mathfrak{m}^{n-1}/\mathfrak{m}^n$ provided $\mathfrak{m}E_n=E_{n-1}$. Let $\mathfrak{m}=(t)$. Then $\mathfrak{m}^{n-1}/\mathfrak{m}^n$ is a 1-dimensional $k$-vector space, $k=R/\mathfrak{m}$, since $\mathfrak{m}^{n-1}=(t^{n-1})$ is principal. The map $\mu_t:E_n \rightarrow E_{n-1}$ has kernel $E_1$. Hence, the image of $\mu_t$ is isomorphic to $N=E_n/E_1$, which is a nonzero submodule of $E_{n-1}$. The $E_n$ are the only submodules of $E$ and $N$ contains elements annihilated by $t^{n-1}$ but no smaller power of $t$. Hence, $N=E_{n-1}=\mathfrak{m}E_n$. Thus,   $E_n/E_{n-1}=E_n/\mathfrak{m}E_n \cong \mathfrak{m}^{n-1}/\mathfrak{m}^n$ is a 1-dimensional $k$-vector space, and $E_n$ is cyclic \cite[Thm 2.3]{Mats}. Note that  if $x_n$ is a generator of $E_n$, then $t^nx_n=0$, but $t^{n-1}x_n \neq 0$. Consequently, $\mathrm{ann}_R(x_n) = \mathfrak{m}^n$ and $E_n \cong R/\mathfrak{m}^n$.

Since $\mu_t:E_n \rightarrow E_{n-1}$ is a surjection for all $n$, there is a special relationship between the generators of the various $E_n$. Let $x_1$ be any generator of $E_1$. We can choose a generator $x_2$ for $E_2$ such that $tx_2=x_1$. Continuing in this fashion, we produce generators $x_1, x_2, x_3, \ldots $ which satisfy
$$tx_1=0, tx_2=x_1, tx_3=x_2, \ldots$$
so that $E$ is a {\em quasicyclic} module (generalizing the notion of quasicyclic group as in \cite[\S3, p 15]{fuchs}), a direct analogue of the group $\mathbb{Z}(p^{\infty}).$

For $R,\mathfrak{m}$ a special $PIR$, the only change necessary is to note that $E=E_n$, $n$ the index of nilpotency of $\mathfrak{m}$. In this case $E$ itself is cyclic.

\begin{propn}
For each prime ideal $P$ of $R$, the submodules $E_n$ of $E_R(R/P)$ are cyclic precisely when $R$ is a finite direct sum of Dedekind domains and special $PIR$s. \qed
\end{propn}

\begin{remk}
The fact that $E=E_R(R/\mathfrak{m})$ is cyclic when $R,\mathfrak{m}$ is a special $PIR$ reaffirms that $R$ is a $0$-dimensional local Gorenstein ring \cite[Prop 21.5]{Eisenbud}. However, there are $0$-dimensional local Gorenstein rings which are {\em not} $PIR$s, for example, $k[x,y]/(x^2,y^3)$ \cite[p 530]{Eisenbud}. Observe that if $R$ is a finite direct sum of Dedekind domains and special $PIR$s, then $R$ is a Gorenstein ring (each localization at a maximal ideal is local Gorenstein \cite[Def p 145]{Mats}). \qed \\
\end{remk}

\end{document}